\documentclass[english,a4,12pt]{article}
\usepackage{a4}
\usepackage{amsmath}
\usepackage{latexsym}
\usepackage{amsfonts}
\usepackage{amssymb}
\usepackage{amscd}
\usepackage{theorem}
\usepackage[dvips]{epsfig}
\usepackage[utf8]{inputenc}
\usepackage{verbatim}
\usepackage{dcolumn}
\usepackage{tabularx}
\usepackage{longtable}
\usepackage{graphics}
\usepackage{pstricks}
\usepackage{pst-node}
\usepackage{latexsym}

\usepackage{babel}
\usepackage{euscript}

\newtheorem{defin}{Definition}[section]
\newtheorem{theorem}{Theorem}

\newtheorem{prop}{Proposition}[section]
\newtheorem{lemma}{Lemma}[section]

\newtheorem{corol}{Corollary}[section]
\newtheorem{rema}{Remark}[section]
\newtheorem{exa}{Example}

\newenvironment{rem}{\begin{rema} \rm}{\end{rema}}

\newenvironment{proof}{{\sc proof:~}}{\hfill $\Box$}

\renewcommand{\le}{\leqslant}

\renewcommand{\ge}{\geqslant}

\renewcommand{\subset}{\subseteq}

\newcommand{\lr}{\leftrightarrow}
\newcommand{\B}{\mathcal{B}}
\newcommand{\tU}{\tilde{U}}
\newcommand{\cU}{\check{U}}
\newcommand{\tZ}{\tilde{Z}}
\newcommand{\tW}{\tilde{W}}
\newcommand{\hY}{\hat{Y}}

\newcommand{\bal}{\begin{align*}}
\newcommand{\eal}{\end{align*}}
\newcommand{\beq}{\begin{eqnarray*}}
\newcommand{\eeq}{\end{eqnarray*}}
\newcommand{\bte}{\begin{theorem}}
\newcommand{\ete}{\end{theorem}}
\newcommand{\bl}{\begin{lemma}}
\newcommand{\el}{\end{lemma}}
\newcommand{\bd}{\begin{description}}
\newcommand{\ed}{\end{description}}

\newcommand{\bc}{\begin{cases}}
\newcommand{\ec}{\end{cases}}
\newcommand{\bp}{\begin{proof}}
\newcommand{\ep}{\end{proof}}
\newcommand{\bco}{\begin{corol}}
\newcommand{\eco}{\end{corol}}

\newcommand{\ue}{\underline{e}}
\newcommand{\ove}{\overline{e}}

\newcommand{\1}{\hbox{1 \hskip -7pt I}}

\newcommand{\iy}{\infty}
\newcommand{\tx}{\text}
\newcommand{\R}{\ensuremath{\mathbb{R}}}

\newcommand{\ov}{\overline}

\newcommand{\Z}{\ensuremath{\mathbb{Z}}}
\newcommand{\N}{\ensuremath{\mathbb{N}}}

\newcommand{\Cst}{\mathsf{Cst}}

\newcommand{\tX}{\tilde{X}}
\newcommand{\hU}{\hat{U}}

\newcommand{\F}{\mathcal{F}}

\newcommand{\Y}{\mathcal{Y}}

\newcommand{\Es}{\mathbb{E}}
\newcommand{\Pb}{\mathbb{P}}
\newcommand{\Var}{\mathsf{Var}}
\newcommand{\Ff}{\mathbb{F}}

\newcommand{\Pc}{\mathcal{P}}

\newcommand{\E}{\mathcal{E}}

\newcommand{\Nn}{\mathcal{N}}

\newcommand{\A}{\mathcal{A}}
\newcommand{\Lc}{\mathcal{L}}

\newcommand{\Gg}{\mathcal{G}}
\newcommand{\M}{{\mathcal M}}
\newcommand{\Qc}{\mathcal{Q}}

\newcommand{\las}{\mathsf{Sp}}

\newcommand{\al}{\alpha}

\newcommand{\be}{\beta}

\newcommand{\g}{\gamma}

\newcommand{\G}{\Gamma}

\newcommand{\De}{\Delta}
\newcommand{\de}{\delta}
\newcommand{\e}{\epsilon}
\newcommand{\Ep}{\Upsilon}

\newcommand{\la}{\lambda}

\newcommand{\s}{\sigma}

\newcommand{\Om}{\Omega}

\def\bdes{\begin{description}}
\def\edes{\end{description}}

\def\iti{\item[(i)]}
\def\itii{\item[(ii)]}

\def\des{\mathop{\hbox{$\longrightarrow$}}\limits}
\def\ses{\mathop{\hbox{$\sim$}}\limits}

\begin{document}         
         
\title{Localization of reinforced random walks}

\author{Pierre Tarr\`es$^1$}         
\maketitle         
 \footnotetext[1]{CNRS, Universit\'e de Toulouse, Institut de Math\'ematiques,
118 route de Narbonne, 31062 Toulouse Cedex 9, France. On leave from the Mathematical Institute, University of Oxford.
E-mail: {\tt tarres@math.univ-toulouse.fr} }
     
\begin{abstract} 
We describe and analyze how reinforced random walks can eventually localize, i.e.  only visit finitely many sites. After introducing vertex and edge self-interacting walks on a discrete graph in a general setting, and stating the main results and conjectures so far on the topic, we present martingale techniques that provide an alternative proof of the a.s. localization of vertex-reinforced random walks (VRRWs) on the integers on finitely many sites and, with positive probability, on five consecutive sites, initially proved by Pemantle and Volkov  (1999,\cite{pemantle2}). 

Next we  introduce the continuous time-lines representation (sometimes called Rubin construction) and its martingale counterpart, and explain how it has been used to prove localization of some reinforced walks on one attracting  edge. Then we show how a modified version of this construction enables one to propose a new short proof of the a.s. localization of VRRWs on five sites on $\Z$.
\end{abstract}

\vspace{0.2in}         
         
{\em AMS 2000 subject classifications.}         
60G50, 60J10, 60K35          
         
{\em Key words and phrases.}         
reinforced walk, martingale, localization

{\em Running Head:}         
Localization of reinforced random walks

\section{Introduction}         
\label{Intro}
Exploration of an environment, behaviour learning or cooperative interaction are instances of situations where the evolution depends on the whole history, either as a tendency to visit again “places” visited before or as a tendency to avoid them. In exploring an unknown city, streets that have been walked before may be considered as more attractive (safer, for instance) or repulsive (boring); learning the best choice among strategies giving random payoffs can be achieved by making random choices with an increasing preference towards the choices that pay more; and cooperation between micro-organisms, for instance, involves miming previously held behaviours. These situations naturally lend themselves to a modelization by self-interacting random processes. 

The definition assumes we are given 
\begin{itemize}
\item $(\Om,\F,\Pb)$ probability space,

\item $(G;\sim)$ nonoriented locally finite graph,

\item $(a_{i,j})_{i,j\in G, i\sim j}$  propensity matrix with positive entries, such that \\
$a_{i,j}>0 \iff i\sim j$,

\item $W:~\N_0\longrightarrow \R_+^*$  weight function.
\end{itemize}

The random process, called $(X_n)_{n\in\N}$, takes values in the set of vertices of $G$; we let $\F_n=\sigma(X_0,\ldots,X_n)$ be the filtration of its past. For all $v\in G$, $n\in\N\cup\{\iy\}$, let 
\begin{equation}
Z_n(v)=\sum_{k=0}^n\1_{\{X_k=v\}}+1
\end{equation}
be the number of visits to $v$ up to time $n$ plus one.

Then $(X_n)_{n\in\N}$ is a  Vertex Self-Interacting Random Walk (VSIRW)
with starting point $v_0\in G$,  propensity matrix $(a_{i,j})_{i\sim j}$ and weight function $W$ 
if 
$X_0=v_0$ and,  for all $n\in\N$, if  $X_n=i$  then
\begin{equation}
\label{def-vsirw}
\Pb(X_{n+1}=j~|~\F_n)=\1_{i\sim j}\frac{a_{i,j}W(Z_n(j))}
{\sum_{k\sim i} a_{i,k}W(Z_n(k))}.
\end{equation}

An Edge Self-Interacting Random Walk (ESIRW) is defined similarly, replacing in \eqref{def-vsirw}  the numbers of visits to  vertices $l\sim i$  by those to the corresponding nonoriented edges $\{i,l\}$: 
\begin{equation}
Z_n(\{i,l\}):=\sum_{k=1}^n(\1_{\{X_{k-1}=i,X_k=l\}}+\1_{\{X_{k-1}=l,X_k=i\}})+1.
\end{equation}

We will define the Edge (resp. Vertex) Reinforced Random Walk as an ESIRW (resp. VSIRW) with linear $W(n)=n+\De$, $\De>-1$: these processes were introduced by Coppersmith and Diaconis in 1986 \cite{coppersmith}.

 In general, the asymptotic behaviour of self-interacting random walks greatly depends on the nature of the interaction. We will focus here on localization phenomena: the difficulty in their analysis lies in the fact that, before this localization occurs, the walk can concentrate on several disconnected clusters -separated by seldom visited sites- so that the relative numbers of visits follow a rather erratic dynamics, which is difficult to analyse. 
 
 For the study of strongly edge reinforced walks (i.e. ESIRW with reciprocally summable weight function $W$), this technical difficulty can be partially overcome by a simple argument, which allows to restrict the study to loop graphs (see Section \ref{rubin}). This argument cannot translate to vertex-reinforced random walks, which display localization on ``richer'' subsets, even on $\Z$, as we describe next.
 \subsection{Preliminary remarks}
 Let us start our study by the following simple preliminary results, which will enable us to gain more intuition on the behaviour of these walks: on one hand on the ``simple'' case of VSIRW on three vertices, and on the other hand on an easy but important property of ESIRWs.
 
 We need to define the two following subsets $R$ and $R'$ of the graph, respectively called range and asymptotic range of the process $(X_n)_{n\ge0}$:
 \begin{align*}
 R&:=\{v\in G\tx{ s.t. } Z_\iy(v)\ne1\}\\
 R'&:=\{v\in G\tx{ s.t. } Z_\iy(v)=\iy\}.
 \end{align*}
 We let $\Cst(x_1,\ldots,x_n)$ be a constant dependent only on $x_1$, $\ldots$, $x_n$. The equalities and inclusions of probability events are understood to hold almost surely.
 \subsubsection{VSIRW on the three consecutive vertices $-1$, $0$ and $1$}  
 \label{3v}
 This walk is equivalent to the ESIRW on two non-oriented edges linking the same pair of vertices, which in turn can be seen as  a $W$-urn process with two colours $-1$ and $1$, defined as follows: we start with a certain number of balls of each colour ($1$ if $X_0=0$) and, at each time step, we pick a ball of colour $i\in\{-1,1\}$ in the urn with a probability proportional to $W(\tx{number of balls of colour }i)$, and put it back  together with a ball of the same colour.
 
 Assume $X_0=0$ for simplicity. Let, for all $n\in\N$, 
 \begin{align}
 \tW_\la(n)&:=\prod_{k=1}^{n-1}\left(1+\frac{\la}{W(k)}\right),\\
 W^*(n)&:=\sum_{k=1}^{n-1}\frac{1}{W(k)},\,\,\hat{W}(n):=\sum_{k=1}^{n-1}\frac{1}{W(k)^2},
 \end{align} 
with the convention that $W^*(1)=\hat{W}(1):=0$ and $\tW_\la(1):=1$. 

For all $\la>0$, let $$A_n(\la):=\frac{\tW_{a_{0,-1}}(Z_n(1))}{\tW_{a_{0,1}}(Z_n(-1))},$$ and let
$$M_n:=\left(\frac{d}{d\la}A_n(\la)\right)_{\la=0}=a_{0,-1}W^*(Z_n(1))-a_{0,1}W^*(Z_n(-1)).$$
\bl
\label{expmart}
The three processes $(A_n(\la))_{n\ge0}$,  $(A_n(\la)^{-1})_{n\ge0}$ and $(M_n)_{n\ge0}$ are martingales.
\el
\bp
Indeed, if $X_n=0$, then 
\begin{align*}
\frac{\Es(A_{n+1}(\la)-A_n(\la)|\F_n)}{A_n(\la)}=&\frac{a_{0,1}W(Z_n(1))}{a_{0,1}W(Z_n(1))+a_{0,-1}W(Z_n(-1))}\frac{a_{0,-1}}{W(Z_n(1))}\\
&-\frac{a_{0,-1}W(Z_n(-1))}{a_{0,1}W(Z_n(1))+a_{0,-1}W(Z_n(-1))}\frac{a_{0,1}}{W(Z_n(-1))}=0.
\end{align*}
\ep

Let us now make use of Lemma \ref{expmart} to analyze the asymptotic behaviour of the walk under the condition $W(n):=(\De+n)^\rho$, $\De>-1$, $\rho\in\R$:
\begin{itemize}
\item[{\bf 1)}] $\rho=1$, i.e. $W(n)=n+\De$.

Then $W^*(n)=\sum_{k=1}^{n-1}1/(k+\De)$, and $W^*(n)-\log n$ converges. On the other hand, for all $n\ge0$,  $$\sum_{k=0}^{n-1}(M_{k+1}-M_k)^2=a_{0,-1}^2\hat{W}(Z_k(1))+a_{0,1}^2\hat{W}(Z_k(-1)),$$ so that 
$(M_n)_{n\ge0}$ is bounded in $L^2$, hence converges a.s. and in $L^2$ by Doob Lemma. Therefore $a_{0,-1}\log(Z_n(1))-a_{0,1}\log(Z_n(-1))$ converges, and 
\begin{equation}
\label{asymp}
Z_n(1)^{a_{0,-1}}\ses_{n\to\iy}C Z_n(-1)^{a_{0,1}}\tx{ a.s.,}
\end{equation}
for a certain positive random variable $C$.

{\bf a)} If $a_{0,1}=a_{0,-1}$, then the $W$-urn is a P\'olya urn, and $\frac{Z_n(1)}{Z_n(1)+Z_n(-1)}$ converges to a random variable $\be\in(0,1)$. 

We can deduce from a classical result that $\be$ is a beta distribution with parameters $(1+\1_{\{X_0=1\}}\De,1+\1_{\{X_0=-1\}}\De)$ in general (but we assumed $X_0=0$ here for simplicity).

{\bf b)} If $a_{0,1}\not=a_{0,-1}$, note that the following urn process with the same colours  $-1$ and $1$ could be analyzed by the technique above: at each time step, pick a ball in the urn with a probability proportional to the number of balls of that colour (as in P\'olya urn), and we put it back together with,  in conditional expectation,  $a_{0,1}$ (resp.$a_{0,1}$) if we picked colour $1$ (resp. $-1$). 

The latter model is similar to Friedman urn \cite{friedman},  analyzed by Freedman in \cite{freedman}. Note that we obtain the same martingale $M_n$ (but not $A_n(\la)$ for general $\la$), and therefore the same asymptotics \eqref{asymp}, if we assume for instance a bounded number of added balls. 

\item[{\bf 2)}] $\rho>1$, or more generally for any $W$ satisfying $\sum1/W(n)<\iy$.

Then the martingale $(M_n)_{n\ge0}$ converges a.s., as a difference of nondecreasing bounded sequences. On the other hand, 
$$\{Z_\iy(1)=Z_\iy(-1)=\iy\}\subset\{M_\iy=0\}.$$

It is possible to prove, using estimates of the variance of the increments, that $\Pb(M_\iy=0)=0$ (see for instance \cite{limic-tarres2}), so that only one of the two vertices $1$ and $-1$ is visited infinitely often almost surely. We will show that result by another technique in Section \ref{rubin}, in Proposition \ref{locurn}.

\item[{\bf 3)}] $\rho<1$ (not necessarily nonnegative).

Given $\la>0$, using  that $(A_n(\la))_{n\ge0}$ is a martingale of expectation $1$, and that, for all $x\ge0$, $0\ge\log(1+x)-x\ge-x^2/2$, 
$$\Es[\exp(\la M_n)]\le\exp(\la^2\hat{W}(n)/2)$$
so that, by Chebychev inequality, for all $a>0$, 
$$\Pb\left[M_n\ge a\sqrt{\hat{W}(n)}\right]\le\exp(-\la a+\la^2\hat{W}(n)/2)\le\exp(-a^2/2),$$ 
choosing $\la:=a/\sqrt{\hat{W}(n)}$. 

Therefore, for all $c>\sqrt{2}$, $|M_n|\le c\sqrt{\log n}\,\hat{W}(n)$ for large $n$ a.s., so that 
$$\frac{Z_n(1)}{Z_n(-1)}\des_{n\to\iy}\left(\frac{a_{0,1}}{a_{0,-1}}\right)^{(1-\rho)^{-1}}$$
\end{itemize}

Hence, on three vertices, the weakly reinforced walk ($\rho\in0,1)$) behaves similarly as the self-repelling one ($\rho<0$), whereas the strongly reinforced walk -i.e. with $W$ reciprocally summable-  implies localization on two vertices. 

In general, is strong reinforcement a necessary and sufficient condition for localization? The next subsection \ref{01} and Section \ref{loce} will  provide a partially positive answer for edge self-interaction (ESIRW), when $W$ is nondecreasing. On the contrary, the results of Section \ref{locv} will highlight the dependence of the behaviour of VSIRWs on the graph structure, which also display localization on particular trapping patterns in the  linear case $W(n)=n+\De$, $\De>-1$. 
\subsubsection{ESIRW on $G$ connected, $W$ nondecreasing, $\sum_{n\in\N}\frac{1}{W(n)}=\iy$}
\label{01}
\begin{prop}
\label{dif}
$\{|R'|\ne0\}=\{R'=G\}$ a.s.
\end{prop}
\bp
Let $t_n:=t_n(x)$ be the $n$-th visit time to $x$, then
$$\sum_{z\sim x} Z_{t_n(x)}(\{x,z\})=2n+a,$$
where $a:=\1_{\{X_0\ne x\}}-2+|\{z\in G: z\sim x\}|$. Hence, for all $z\sim x$ and $n\in\N$,
\begin{align*}
&\Pb(X_{t_n+1}=z|\F_{t_n})\1_{\{t_n<\iy\}}\ge\frac{a_{x,z}}{\sum_{y\sim x}a_{x,y}}\frac{W(0)}{W(2n+a)}\\
&\ge \1_{\{t_n<\iy\}}W(0)\Cst(a_{x,y})_{y\sim x}
\left(\frac{1}{W(2n+a)}+\frac{1}{W(2n+1+a)}\right),
\end{align*}
using that $W$ is nondecreasing.

Therefore, using conditional Borel-Cantelli Lemma \ref{bc}, 
\begin{align*}
\{Z_\iy(x)=\iy\}\subset&\left\{\sum_{n\in\N}\Pb(X_{t_n+1}=z|\F_{t_n})\1_{\{t_n<\iy\}}=\iy\right\}\\
&=\left\{\sum_{n\in\N}\1_{\{X_{t_n+1}=z\}}\1_{\{t_n<\iy\}}=\iy\right\}\subset\{Z_\iy(z)=\iy\}\tx{ a.s.}
\end{align*}
\ep
\begin{rem}
The condition that $W$ should be nondecreasing is important in Proposition \ref{dif}. The following counterexample was proposed by Sellke \cite{sellke2}: if $\sum\frac{1}{W(2k)}=\iy$ and $\sum\frac{1}{W(2k+1)}<\iy$, $a_{i,j}=\1_{\{i\sim j\}}$, $G=\R^d$ and $X_0=0$, then 
$$\Pb(\forall n\in\N, \,X_{2n}=0)>0.$$
This result will be a direct consequence of the time-lines construction in Section \ref{rubin}.
\end{rem}
\subsection{Localization results, VSRIW}
\label{locv}
\subsubsection{$G=\Z$, $a_{i,j}=\1_{i\sim j}$, $W(n)=\De+n$}
Pemantle and Volkov \cite{pemantle2} showed that the walk a.s. visits only finitely many vertices, and that it localizes on any set of five consecutive sites with positive probability. Tarrès \cite{tarres2} proved that localization on fives sites is in fact the a.s. behaviour.
\bte [Pemantle and Volkov, \cite{pemantle2}]
\label{loc} 
$|R'|<\iy$ a.s. and, for any $x\in\Z$, $\Pb(R'=\{x-2,x-1,x,x+1,x+2\})>0$.
\ete
\bte [Tarrès, \cite{tarres2}]
\label{loc5}
$|R'|=5$ a.s. 
\ete 
We propose here some new proofs of these results: Theorem \ref{loc} in Section \ref{1stres} (see Propositions \ref{pemvolk} and \ref{bienvenue}), using techniques derived from \cite{tarres2}, and Theorem \ref{loc5} in Section \ref{short}, partially adapting a continuous-time equivalent of the random walk, originally introduced in \cite{davis,sellke2}.
\subsubsection{$G=\Z$, $a_{i,j}=\1_{i\sim j}$, $W(n)\sim n^\al$}
\bte [Volkov, \cite{volkov2}]
Suppose that $0<\lim_{k\to\iy} w_k/k^{\al}<\iy$. Then
\begin{description}
\item[(a)] If $\al<1$, then $R'\in\{0,\iy\}$
\item[(b)] If $\al>1$, then $|R'|=2$.
\end{description}
\ete
\subsubsection{General $(G,\sim)$ locally finite, $(a_{i,j})_{i\sim j}$ symmetric and $W(n)=\De+n$}
Under a symmetric propensity matrix, the vertex-reinforced random walk localizes with positive probability on a class of complete $d$-partite subgraphs with possible loops plus their outer boundary. We need to introduce some notation, in order to describe further these trapping subsets.

Given a subset $R$ of $G$, we let 
$$\partial R=\{j\in G\setminus R~:~j\sim R\}$$
be  the {\it outer boundary} of $R$.

For any $x=(x_i)_{i\in G}\in\R^{G}$, let 
$$S(x):=\{i\in G/~x_i\not=0\}$$ be its support. Let 
$$\De:=\left\{x\in\R_+^{G}\tx{ s.t. }|S(x)|<\iy\tx{ and }\sum_{i\in G}x_i=1\right\}$$
be the nonnegative simplex restricted to elements $x$ of finite support.

For all $x\in\De$, let 
\begin{equation}
\label{nih}
N_i(x):=\sum_{j\in G, j\sim i}a_{i,j}x_j,~~~~
H(x)=\sum_{i,j\in G, i\sim j}a_{i,j}x_ix_j=\sum_{i\in G}x_iN_i(x).
\end{equation}
For all $n\in\N$, let 
$$x_n=\left(\frac{Z_n(i)-1}{n}\right)_{i\in G}$$
be the vector of density of occupation  of the random walk at time $n$, which has finite support and takes values in $\De$.

The following definition introduces ``good candidates'' for the limiting density of occupation of the random walk.
\begin{defin}
For all $x\in\De$, let ${\bf(P)}_{x}$ be the following predicate:
$$\max \left (\las\left[a_{i,j}-2H(x)\right]_{i,j\in S(x)}\right)\le0,\,\,\max\{N_i(x)-H(x),i\in\partial S(x)\} < 0.$$
\bte [Benaïm and Tarrès, \cite{benaimtarres}]
Given $x\in\De$ such that such that ${\bf(P)}_{x}$ holds and for any neighbourhood $\Nn(x)$ of $x$ in $\De$, there is with positive probability $y\in\Nn(x)$ with $S(y)=S(x)$, such that the following three events occur:
\bal
{\bf(i)}&\,\,\,R'=S(x)\cup\partial S(x)\\
{\bf(ii)}&\,\,\,x_n\to y\\
{\bf(iii)}&\,\,\,\forall i\in\partial S(x),~Z_n(i)/n^{N_i(x)/H(x)}\to C_i\in (0,\iy)\tx{ (random)}.
\end{align*}
\ete
\end{defin}
Assumption ${\bf(P)}_{x}$ in Definition \ref{nih} describes stable equilibria of the ordinary differential equation
\begin{equation}
\label{eqdiff}
\frac{dx}{dt}=F(x),
\end{equation}
where 
\begin{equation}
\label{replicator}
F(x)=(x_i[N_i(x)-H(x)])_{i\in G}.
\end{equation}
also known as the linear {\it replicator }equation in population genetics and game theory. 
Up to an adequate rescaling in time, we can indeed show that $(x_k)_{k\in\N}$ approximate of this ODE under certain assumptions. 

The support of these equilibria satisfies some properties, described in the following Lemma \ref{descr}. In the context of population dynamics, they inform on the structure of the surviving species, depending on the nature of the graph.
\begin{defin}
Given $d\ge1$, subgraph $(S,\sim)$ of $(G,\sim)$ will be called a {\em complete $d$-partite graph with possible loops}, if $(S,\sim)$ is a $d$-partite graph on which some loops have possibly been added. That is
 $$S=V_1\cup\ldots\cup V_d$$ with 
\bdes
\iti
 $\forall$ $p\in\{1,\ldots,d\}$, $\forall$ $i$, $j$ $\in V_p$, if $i\not=j$ then $i\not\sim j.$

\itii  $\forall$ $p$, $q$ $\in\{1,\ldots,d\}$, $p\not=q$, $\forall i\in V_p$, $\forall j\in V_q$, $i\sim j$. 
\edes
\end{defin}
\begin{defin}
For all $S\subset G$, let ${\bf(P)}_S$ be the following predicate:
\bal
&{\textnormal{\bf(P)}_S{\bf(a)}}~\tx{
 $(S,\sim)$ is a complete $d$-partite graph with possible loops. }\\
&{\textnormal{\bf(P)}_S{\bf(b)}}~\tx{If $i\sim i$ for some $i\in S$, then the partition containing $i$ is a singleton.}\\
&{\textnormal{\bf(P)}_S{\bf(c)}}~\tx{If $V_p$, $1\le p\le d$ are its $d$ partitions, then for all 
$p$, $q$ $\in\{1,\ldots,d\}$ and }\\
&~~~~~~~~~~~~~~~~\tx{$i$, $i'$ $\in V_p$, $j$, $j'$ $\in V_j$, 
$a_{i,j}=a_{i',j'}$. }
\end{align*}
\end{defin}
\bl
\label{descr} (Benaïm and Tarrès, \cite{benaimtarres})
For all $x\in\De$, ${\bf(P)}_{x}$ implies ${\bf(P)}_{S(x)}$.
\el
\subsection{Localization results, ESIRW}
\label{loce}
Let (H) be the following condition on $W$:         
$$\sum_{k\in\N}\frac{1}{W(k)}<\iy.$$
We know from Proposition \ref{dif} that, if (H) does not hold and $W$ is nondecreasing, then the walk is either transient or recurrent on all vertices, assuming the graph $G$ is connected.          
      
 On the other hand, it is easy to show that (H) implies localization on a single edge with positive probability. Sellke \cite{sellke2} conjectured in 1994 that this should occur with probability one on any graph of bounded degree, and proved the statement on graphs without odd cycles. Then Limic and Tarrès \cite{limic-tarres}  showed in 2007 that this conjecture indeed holds if $W$ is nondecreasing (Limic \cite{limic} solved the case $W(k)=(k+1)^\rho$ in 2003).
 
 \begin{theorem} [Sellke \cite{sellke2}, Theorem 3, Limic \cite{limic}, Lemmas 1-2, Corollaries 1-2]
\label{sellke}
If $(G,\sim)$ has bounded degree and contains no odd cycles, then (H) implies $|R'|=2$ a.s.
\end{theorem}

 \begin{theorem} [Limic and Tarrès \cite{limic-tarres}, Corollary 3]
\label{limtar}
If $(G,\sim)$ has bounded degree and
$W$ is nondecreasing, then (H) implies $|R'|=2$ a.s.
\end{theorem}
We explain the key ingredient for the proof of Theorem \ref{sellke} in Section \ref{erubin}.
 \section{VRRW on $\Z$, first localization results}
 \label{1stres}
 The goal of this section is to prove the following Propositions \ref{pemvolk} and \ref{bienvenue}, which will in particular imply Theorem \ref{loc}.  \begin{prop} (Pemantle and Volkov, 1999, \cite{pemantle2})
 \label{pemvolk}
 \begin{align*}
 {\textnormal{\bf{(a)}}}&~\tx{ For all $x\in\Z$, $\al\in(0,1)$, $\e>0$,}\\
&\,\,\,\,\, \Pb(\{R'=\{x-2,x-1,x,x+1,x+2\}\}\cap\{\al_\iy^-(x)\in(\al-\e,\al+\e)\})>0.\\
 {\textnormal{\bf{(b)}}}&~|R|<\iy\tx{ a.s. }
 \end{align*}
 \end{prop}
 \begin{prop}
 \label{bienvenue}
 \bal
  {\textnormal{\bf{(a)}}}&~ \tx{(Bienvenüe, 1999, \cite{bienvenue}) }\\
  &\{R'=\{x-2,x-1,x,x+1,x+2\}\}\\
  &\subset\{\exists\al_\iy^-(x):=\lim\al_n^-(x)\in(0,1)\}\cap\{\log Z_n(x\pm2)/\log Z_n(x)\des_{n\to\iy}\al_\iy^\pm(x)\}\tx{ a.s.}\\
   {\textnormal{\bf{(b)}}}&~ \tx{(Pemantle and Volkov, 1999, \cite{pemantle2})}
   |R'|\ge5\tx{ a.s. }
  \end{align*} 
 \end{prop}
 
 We will assume that $\De:=0$ in the remainder of this survey, so that $W(n)=n$; the proofs obviously carry on to the general case, by replacing $Z_.(.)$ by $Z_.(.)+\De$, and by defining the following function $h$ accordingly, on $\De+\N$ instead of $\N$. Also recall that $a_{i,j}=\1_{i\sim j}$ here.

Let us first explain the heuristics of the localization on five vertices with positive probability. Let us assume for instance that the that the walk started at $0$, and that we are in the following configuration at time $n$, represented by the figure below: site $0$ has been visited $n/2$ times, its neighbours $-1$ and $1$ have shared the other half of time, and have been visited respectively roughly $\al n/2$ and $(1-\al)n/2$ times for some constant $\al\in(0,1)$; that the two sites $-2$ and $2$ have been visited of the order of $Cn^\al$ and $C'n^{1-\al}$ and, finally that $-3$ and $3$ have not been visited yet. 

The numbers above the sites on the figure represent the estimates of the numbers of visits (plus one).

\bigskip
\bigskip
\bigskip
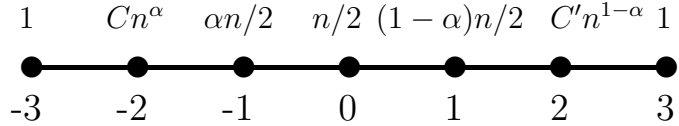
\begin{picture}(400,40)
\linethickness{1.5pt}
\put(60,30){\line(1,0){240}}
\put(60,30){\circle*{8}}
\put(100,30){\circle*{8}}
\put(140,30){\circle*{8}}
\put(180,30){\circle*{8}}
\put(220,30){\circle*{8}}
\put(260,30){\circle*{8}}
\put(300,30){\circle*{8}}
\put(52,10){\large{-3}}
\put(92,10){\large{-2}}
\put(132,10){\large{-1}}
\put(176,10){\large{0}}
\put(216,10){\large{1}}
\put(256,10){\large{2}}
\put(296,10){\large{3}}
\put(55,45){$1$}
\put(88,45){$Cn^\al$}
\put(125,45){$\al n/2$}
\put(166,45){$n/2$}
\put(190,45){$(1-\al) n/2$}
\put(256,45){$C'n^{1-\al}$}
\put(296,45){$1$}
\end{picture}

Assume $X_n=-1$ for instance. First, the configuration would not be viable asymptotically for an edge-reinforced random walk (ERRW), since 
$$\Pb(X_{n+2}=-3|\F_n)\ses_{n\to\iy}\frac{2Cn^\al}{\al n}\frac{1}{2Cn^\al}\ses_{n\to\iy}\frac{1}{\al n}$$
so that $-3$ would eventually be visited, by conditional Borel-Cantelli lemma. 

On the contrary, for the vertex-reinforced random walk we are considering, the same computation yields
$$\Pb(X_{n+2}=-3|\F_n)\ses_{n\to\iy}\frac{Cn^\al}{n/2}\frac{1}{\al n/2}\ses_{n\to\iy}\frac{\Cst(C,\al)}{n^{2-\al}}.$$

Therefore, with lower bounded probability, the sites $-3$ and $3$ will never be visited, as long as the same asymptotics holds for the visits to the other sites. Now, under these assumptions:
\begin{itemize}
\item The visits to $-2$ and $2$ are seldom, so that the respective visits to $-1$ and $1$ can almost be estimated by considering the walk on the three vertices $-1$, $0$ and $1$, restricted to its moves to and from $0$, described in Section \ref{3v}: we are in Case 1)a), so that $Z_n(-1)/(Z_n(-1)+Z_n(1))$ should remain close to $\al$.
\item The respective visits to $-2$ (similarly $2$) and $0$ can be seen as stemming from a Friedman urn model, Case 2)a) of Section \ref{3v}. Indeed, starting from $-1$,  the sites $-2$ and $0$ are chosen proportionally to their numbers of visits. But, if $-2$ is picked then the walk immediately comes back to $-1$, whereas if $0$ is chosen, the expectation of the number visits before coming back to $-1$ is roughly $\al^{-1}$: \eqref{asymp} would provide convergence of $Z_n(-2)/Z_n(0)^\al$. 
\end{itemize}

The following results indeed justify this heuristics, and will also be useful in the proof of a.s. localization on five points. They are similar to those developped in Sections 1, 2 and 3 of \cite{tarres2}.

  Let, for all $n\in\N$, 
  $$h(n):=\sum_{k=1}^{n-1}\frac{1}{k},$$
 with the convention that $h(1):=0$. For all $n\in\N_0$ and $x\in\Z$, denote
\begin{align*}
Z_n^\pm(x)&:=\sum_{k=1}^{n} \1_{\{X_{k-1}=x,X_k=x\pm1\}},\\
\al_n^\pm(x)&:=\frac{Z_n(x\pm1)}{Z_n(x-1)+Z_n(x+1)},\\
Y_n^\pm(x)&:=\sum_{k=1}^n \1_{\{X_{k-1}=x,X_k=x\pm1\}}
\frac{1}{Z_{k-1}(x\pm1)},\\
\end{align*}
also
\begin{align*}
Y_n(x)&:=\sum_{k=1}^n \1_{\{X_{k-1}=x\}}\frac{1}{Z_{k-1}(x-1)+Z_{k-1}(x+1)},\\
\hY_n^\pm(x)&:=Y_n^\pm(x)-Y_n(x),
\end{align*}
which are respectively the previsible and martingale part in the Doob decomposition of $Y_n^\pm(x)$, 
and finally
$$Y_\iy^\pm(x):=\lim_{n\to\iy}Y_n^\pm(x),\,\,Y_\iy(x):=\lim_{n\to\iy}Y_n(x).$$

Given $(a_n)$, $(b_n)$ random processes on $\R$, we write $a_n\equiv b_n$ iff $a_n-b_n$ converges a.s.

Let us define the probability event
$$\Ep(x):=\{Y_\iy(x)<\iy\}$$
and, for any finite sequence 
$(x_i)_{1\le i\le n}$ taking values in $\Z$, 
the event
$$\Ep( (x_i)_{1\le i\le n})=\bigcap_{1\le i\le n}\Ep(x_i).$$

\subsection{``P\'olya urn'' estimates}
The event $\Ep(x)$ corresponds to the event that $x$ is ``seldom'' visited (represented by a cross on the figure), hence ``neutral'' with respect to its neighbours, in the following sense: the respective visits to $x+1$ and $x+3$ starting from $x+2$ can be seen be seen similar to those of a P\'olya  urn model (see Section \ref{3v}, study of case 1)a)), perturbed by the visits from $x$ and $x+4$: $\Ep(x)$ implies that the visits from $x$ do not act upon the asymptotic behaviour of $\al_n^-(x+2)$, as stated in Corollary \ref{noinf}.
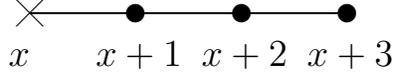
\begin{center}
\begin{picture}(300,40)
\put(100,30){\line(1,0){120}}
\put(95,35){\line(1,-1){10}}
\put(95,25){\line(1,1){10}}
\put(140,30){\circle*{7}}
\put(180,30){\circle*{7}}
\put(220,30){\circle*{7}}
\put(92,10){\large{$x$}}
\put(125,10){\large{$x+1$}}
\put(165,10){\large{$x+2$}}
\put(205,10){\large{$x+3$}}
\end{picture}
\end{center}
\bigskip
 \begin{prop}
 \label{1est}
 For all $x\in\Z$, 
 \bal
&{\textnormal{\bf{(a)}}}~\hY_n^\pm(x)=Y_n^\pm(x)-Y_n(x)\tx{ is a martingale, converging a.s. and in $L^2$}\\
&{\textnormal{\bf{(b)}}}~Y_n^\pm(x)\equiv Y_n(x)\\
&{\textnormal{\bf{(c)}}}~\Es((\hY_n(x)-\hY_\iy(x))^2|\F_n)\le\Cst\, Z_n(x\pm1)^{-1}\\
&{\textnormal{\bf{(d)}}}~Y_n^+(x-1)+Y_n^-(x+1)=h(Z_n(x))-h(1+{\textnormal{\1}}_{\{X_0=x\}})\equiv\log Z_n(x)
\end{align*}
\end{prop}
\bp
It follows from its definition that $(\hY_n^\pm(x))_{n\ge0}$ is a martingale. Now
\bal
\Var(\hY_{n+1}^\pm(x)|\F_n)=\Var(Y_{n+1}^\pm(x)|\F_n)\le\,\, &\Es((Y_{n+1}^\pm(x)-Y_{n}^\pm(x))^2|\F_n)\\
&=\Es\left(\frac{\1_{\{X_n=x,X_{n+1}=x\pm1\}}}{Z_{n}(x\pm1)^2}|\F_n\right)
\end{align*}
Hence, for all $m\ge n$, 
\bal
\Es((Y_{n}^\pm(x)-Y_{m}^\pm(x))^2|\F_n)\le\Es\left(\sum_{k=n}^\iy\frac{\1_{\{X_k=x,X_{k+1}=x\pm1\}}}{Z_{k}(x\pm1)^2}|\F_n\right)
\le\sum_{l=Z_n(x\pm1)}^\iy\frac{1}{l^2}\le\frac{\Cst}{Z_n(x\pm1)}.
\end{align*}
This implies {\bf (a)-(c)}; {\bf (d)} follows from definitions.
\ep

Let, for all $t\in(0,1)$, 
$$f(t):=\log\left(\frac{t}{1-t}\right).$$
Let, for all $a>0$, $m\ge0$, $y\in\Z$, let $\E_{m,a}(y)$ be the event
\begin{equation}
\label{emay}
\E_{m,a}(y):=\left\{\sup_{m\le k\le n<\iy}[f(\al_n^-(y))-f(\al_k^-(y))-(Y_n(y-2)-Y_k(y-2))]\ge a\right\}.
\end{equation}
\bco
\label{noinf}
 For all $x\in\Z$, $n\ge m$, $a\ge\Cst$, 
 \bal
&{\textnormal{\bf{(a)}}}~\Ep(x)\subset\{\exists\al_\iy^-(x+2):=\lim\al_n^-(x+2)\in[0,1)\}\cap\{\log f(\al_n^-(x+2))\equiv -Y_n^-(x+4)\}\\
&{\textnormal{\bf{(b)}}}~\Ep(x)\cap\{\al_\iy^-(x+2)>0\}\subset\Ep(x+4)\\
&{\textnormal{\bf{(c)}}}~\Pb(\E_{m,a}(x+2)|\F_m)\le\Cst(a)\left(Z_m(x+1)^{-1}+Z_m(x+3)^{-1}\right)\\
&{\textnormal{\bf{(d)}}}~\Ep(x-1,x+1)\subset\{Z_\iy(x)<\iy\}
\end{align*}
\eco
\bp
{\bf(a)-(b)} By Proposition \ref{1est} {\bf(b)} and {\bf(d)}, a.s. on $\Ep(x)=\{Y_\iy^=(x)<\iy\}$, 
\bal
\log Z_n(x+1)\equiv Y_n^+(x)+Y_n^-(x+2)\equiv Y_n^+(x+2)\equiv\log Z_n(x+3)-Y_n^-(x+4),
\end{align*}
so that 
$$\log f(\al_n^-(x+2))=\log\frac{Z_n(x+1)}{Z_n(x+3)}\equiv-Y_n^-(x+4).$$

{\bf(c)} By Proposition \ref{1est} {\bf(d)},
\begin{align}
\nonumber
&h(Z_n(x+1))-h(1+\1_{\{X_0=x+1\}})
=Y_n^+(x)+Y_n^-(x+2)\\
\nonumber
&=Y_n(x)+Y_n^+(x+2)+\hY_n^+(x)+\hY_n^-(x+2)-\hY_n^+(x+2)\\
\nonumber
&=Y_n(x)+h(Z_n(x+3))-h(1+\1_{\{X_0=x+3\}})-Y_n^-(x+4)\\
\label{longest}
&\,\,\,\,\,\,\,\,\,\,\,\,\,\,\,\,\,\,\,\,\,\,\,\,\,\,\,\,\,\,\,\,\,\,\,\,\,\,\,\,\,\,\,\,\,\,\,\,\,\,\,\,\,\,\,\,\,+\hY_n^+(x)+\hY_n^-(x+2)-\hY_n(x+2).
\end{align}
In order to estimate the probability of $\E_{m,a}(x+2)$ we observe that, if $a\ge\Cst$, then
\bal
&f(\al_n^-(x+2))-f(\al_m^-(x+2))=\log\left(\frac{Z_n(x+1)}{Z_m(x+1)}\right)-\log\left(\frac{Z_n(x+3)}{Z_m(x+3)}\right)\\
&\le h(Z_n(x+1))-h(Z_m(x+1))-[h(Z_n(x+3))-h(Z_m(x+3))]+\frac{a}{4},
\end{align*}
using that, for any large $q\ge p$, 
\begin{equation}
\label{logapprox}
-\frac{\Cst}{p}\le\log\frac{q}{p}-(h(q)-h(p))\le0,
\end{equation}
as $1/p-1/p^{2}\le\log((p+1)/p)\le 1/p$ for large $p$. 

We subtract identities \eqref{longest} at times $n$ and $m$, and note that 
the terms $\hY_n^+(x)-\hY_m^+(x)$, $\hY_n^-(x+2)-\hY_m^-(x+2)$ and $\hY_m(x+2)-\hY_n(x+2)$ can be upper bounded by $a/4$ with large probability, by Doob and Chebychev inequalities, and Proposition \ref{1est} {\bf(c)}, for instance:
\bal
\Pb\left(\sup_{n\ge m}(\hY_n^+(x)-\hY_m^+(x))\ge \frac{a}{4}|\F_m\right)\le\frac{\Cst}{a^2}\Es((\hY_n(x)-\hY_\iy(x))^2|\F_n)\le\frac{\Cst}{a^2}\frac{1}{Z_m(x\pm1)},
\end{align*}
which completes the proof.

{\bf(d)} It is a direct consequence of Proposition \ref{1est} {\bf(b)} and {\bf(d)}: a.s. on $\Ep(x-1,x+1)$, 
$$\log Z_\iy(x)=Y_\iy^+(x-1)+Y_\iy^-(x+1)<\iy.$$
\ep
\subsection{``Friedman urn'' estimates}
The main goal of this subsection is to provide estimates for the number of visits to $x+1$ when $\al_n^-(x+3)$ converges to $\al_\iy^-(x+3)>0$, and $\Ep(x-1)$ holds, i.e. $x-1$ is ``seldom'' visited. The Friedman urn dynamics, which arises for sites $x$ and $x+2$, appears through the following calculation, justified rigourously in this section:
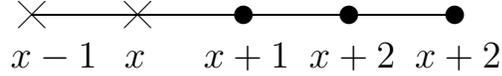
\begin{center}
\begin{picture}(300,40)
\put(100,30){\line(1,0){160}}
\put(95,35){\line(1,-1){10}}
\put(95,25){\line(1,1){10}}
\put(135,35){\line(1,-1){10}}
\put(135,25){\line(1,1){10}}
\put(180,30){\circle*{7}}
\put(220,30){\circle*{7}}
\put(260,30){\circle*{7}}
\put(92,10){\large{$x-1$}}
\put(135,10){\large{$x$}}
\put(165,10){\large{$x+1$}}
\put(205,10){\large{$x+2$}}
\put(245,10){\large{$x+2$}}
\end{picture}
\end{center}
\medskip
\bal
&\log Z_n(x)\equiv Y_n^-(x+1)\equiv Y_n^+(x+1)=\sum_{k=1}^n\frac{\1_{\{X_{k-1}=x+1,X_k=x+2\}}}{Z_{k-1}(x+2)}\\
&\approx
\sum_{k=1}^n\frac{\1_{\{X_{k-1}=x+2,X_k=x+1\}}}{Z_{k-1}(x+2)}\approx
\sum_{k=1}^n\frac{\1_{\{X_k=x+2\}}}{Z_{k-1}(x+2)}\al_k^-(x+2)\approx\al_\iy^-(x+2)\log Z_n(x+2).
\end{align*}

For all $x\in\Z$, $n\in\N$, $\al>0$, let 
\bal
U_{n,\pm}^+(x)&:=\sum_{k=1}^n\frac{\1_{\{X_{k-1}=x,X_k=x\pm1\}}}{Z_{k-1}(x)}\\
U_{n,\pm}^-(x)&:=\sum_{k=1}^n\frac{\1_{\{X_{k-1}=x,X_k=x\mp1\}}}{Z_{k-1}(x)}\frac{Z_{k-1}(x+1)}{Z_{k-1}(x-1)},
\end{align*}
and let 
\bal
U_{n,\pm}(x)&:=\sum_{k=1}^n\frac{\1_{\{X_{k-1}=x\}}}{Z_{k-1}(x)}\al_{k-1}^\pm(x)\\
\hU_{n,\pm}^\pm(x)&:=U_{n,\pm}^\pm(x)-U_{n,\pm}(x)
\end{align*}
be respectively the martingale and previsible parts of the Doob decomposition of $U_{n,\pm}^\pm(x)$. 

Note that the processes $U_{n,\pm}^-(x)$ and $\hU_{n,\pm}^-(x)$ will not be used in this section; they will be necessary tools to the proof of a.s. localization in Section \ref{short}.

\begin{prop}
\label{pol}
For all $x\in\Z$, $C$, $a$, $\g$ $>0$, $n\in\N$, 
\bal
&{\textnormal{\bf{(a)}}}~n\mapsto Y_n^\pm(x)+\frac{{\textnormal{\1}}_{\{\pm X_n\le\pm x\}}}{Z_{n-1}(x\pm1)}-U_{n,\mp}^+(x\pm1)\tx{ is a.s. constant}\\
&{\textnormal{\bf{(b)}}}~\hU_{n,\pm}^+(x)\tx{ converges a.s. and in $L^2$, and}\\
&\,\,\,\,\,\,\,\,\,\,\,\,\,\,\Es((\hU_{n,\pm}^+(x)-\hU_{\iy,\pm}^+(x))^2|\F_n)\le\Cst\, Z_n(x)^{-1}\\
&{\textnormal{\bf{(c)}}}~\Ep(x-1)\cap\{Z_\iy(x)=\iy\}\cap\{\limsup\al_n^-(x+2)\le\g\}\\
&\,\,\,\,\,\,\,\,\,\,\,\,\,\,\subset\Ep(x-1,x)\cap\{\log Z_n(x)\ses_{n\to\iy}\al_\iy^-(x+2)\log Z_n(x+2)\}\\
&{\textnormal{\bf{(d)}}}~\tx{ Assume }Z_m\le CZ_m(x+2)^\al\tx{, and let}\\
&T:=\inf\{n\ge m\tx{ s.t. }\al_n^-(x+2)\ge\g\tx{ or }X_n=x-1\}.\tx{ Then}\\
&\Pb\left(\sup_{T>n\ge m}\frac{Z_n(x)}{Z_n(x+2)^\g}\ge e^{a}C|\F_m\right)\le\Cst(a)(Z_m(x)^{-1}+Z_m(x+2)^{-1})
\end{align*}
\end{prop}
\bp
{\bf(a)} Assume $\pm=+$ for simplicity. Now the function of $n$ only changes when

-  $X_{n-1}=x$ and $X_n=x+1$, in which case $Y_n^+(x)$ increases by $Z_{n-1}(x+1)^{-1}$

-  $X_{n-1}=x+1$ and $X_n=x$, in which case $U_{n,-}^+(x+1)$ increases by $Z_{n-1}(x+1)^{-1}$

{\bf(b)} is similar to Proposition \ref{1est} {\bf(c)}

{\bf(c)} Note that
\bal
\log Z_n(x)&\equiv Y_n^+(x-1)+Y_n^-(x+1)\equiv Y_n^+(x-1)\equiv U_{n,-}^+(x+2)\equiv U_{n,-}(x+2)\\
&\equiv\sum_{k=1}^n\frac{\1_{\{X_{k-1}=x+2\}}}{Z_{k-1}(x+2)}\le(\g+\e)\log Z_n(x+2)
\end{align*}
for large $n$ if $\limsup\al_n^-(x+2)\le\g$ and $Z_\iy(x)=\iy$. 
This implies $\al_k^-(x+1)\le(Z_n(x)+Z_n(x+2))^{\g+\e-1}
\le Z_n(x+1)^{\g+\e-1}$ for large $n\in\N$, so that 
$$Y_\iy(x)\equiv Y_\iy^+(x)\equiv U_{n,-}(x+1)=\sum_{k=1}^n\frac{\1_{\{X_{k-1}=x+1\}}}{Z_{k-1}(x+1)}\al_{k-1}^-(x+1)<\iy.$$

{\bf(d)} Similarly as in Corollary \ref{noinf} {\bf(c)}, outside of an event of probability less than $\Cst(a)(Z_m(x)^{-1}+Z_m(x+2)^{-1})$, we have
\bal
\log\left(\frac{Z_n(x)}{Z_m(x)}\right)\le
Y_n^+(x-1)-Y_m^+(x-1)+a+\sup_{m\le k\le n}\al_{k-1}^\pm(x+2)\log\left(\frac{Z_n(x+2)}{Z_m(x+2)}\right).
\end{align*}
\ep
\subsection{Conclusions}
\subsubsection{Proof of Proposition \ref{bienvenue}}
{\bf(a)} A.s. on $\{Z_\iy(x+3)<\iy\}\cap\{Z_\iy(x-3)<\iy\}$, 
$\Ep(x-3,x-2,x+2,x+3)$ holds by Proposition \ref{1est} {\bf(b)}. Hence $\al_n^-(x)\des_{n\to\iy}\al_\iy^-(x)\in[0,1)$ by Corollary  \ref{noinf} {\bf(a)}, which also implies $\al_\iy^-(x)\in[0,1)$ by considering the sites in the reverse order, so that $\al_\iy^-(x)\in(0,1)$. Finally Proposition  \ref{pol} {\bf(c)} completes the proof.

{\bf(b)} A.s. $\{R'\subset\{x,x+1,x+2,x+3\}\}\cap\{Z_\iy(x)=\iy\}$, Corollary \ref{noinf} {\bf(a)} implies $\al_n^-(x+2)\des_{n\to\iy}\al_\iy^-(x+2)\in[0,1)$; subsequently,  $\al_n^-(x+1)\des_{n\to\iy}0$ by Proposition \ref{pol} {\bf(c)}. On the other hand, $\Ep(x+3)$ holds, so that $\al_\iy^+(x+1)\des_{n\to\iy}\al_\iy^+(x+1)\in[0,1)$, which is contradictory.
\subsubsection{Proof of localization with positive probability}
We will prove the following result of localization on a half-axis.
\bco
\label{cor-loc}
For all $x\in\Z$, $\g\in[0,1)$, $\e>0$, $C>0$, $m\in\N$, assume $Z_m(x\pm3)=1$, $Z_m(x\pm2)\le CZ_m(x)^\g$ and $\al_m^\pm(x)\le\g$. Then
$$\Pb\left(\{Z_\iy(x\pm3)\ne1\}\cup\{\sup_{n\ge m}\al_n^\pm(x)>\g+\e\}|\F_m\right)\le\Cst(\g,\e,C)(Z_m(x)^{-\g}+Z_m(x\pm1)^{\g+\e-1}).$$
\eco
\bp
Assume $\pm:=-$, $X_m\ge x-1$ and $\e<1-\g$ for simplicity. Let $T_1$, $T_2$ and $T_3$ be the stopping times
\bal
T_1&:=\inf\{n\ge m\tx{ s.t. }Z_n(x-3)>1\},\\
T_2&:=\inf\{n\ge m\tx{ s.t. }\al_n^-(x)\ge\g+\e\}\\
T_3&:=\inf\{n\ge m\tx{ s.t. }Z_n(x-2)\ge Ce^\e Z_n(x)^{\g+\e}\}.
\end{align*}
First, if $Z_m(x-2)\ge\Cst(C,\e)$,
\begin{align}
\nonumber
\Pb(T_1<T_2\wedge T_3|\F_m)&\le\sum_{k=m}^\iy\1_{\{X_k=x-1\}}\frac{Ce^\e Z_k(x)^{\g+\e}}{Z_k(x-2)+Z_k(x)}\frac{1}{Z_k(x-1)}\\
\label{t1}
&\le 2Ce^\e\sum_{k=Z_m(x-1)}^\iy k^{\g+\e-2}\le\frac{2Ce^\e(1-(\g+\e))^{-1}}{(Z_m(x-1)-1)^{1-(\g+\e)}}
\end{align}

Second, 
\begin{equation}
\label{t2}
\Pb(T_2<T_1\wedge T_3|\F_m)\le\Pb(\E_{m,\e/2}(x)|\F_m),
\end{equation}
where $\E_{.,.}(.)$ is defined in \eqref{emay}. Indeed, assume $\E_{m,\e/2}(x)$ holds. For all $n\ge m$, let $p(n)$ be the last time $k\le n$ s.t. $\al_k^-(x)\le\g$. Then
\bal
Y_{p(n),n}(x-2)\le\g^{-1}\sum_{k=p(n)+1}^n\frac{\1_{\{X_k=x-2\}}}{Z_k(x)}
\le\g^{-1}(Ce^\e)^{(\g+\e)^{-1}}\sum_{k=p(n)+1}^n\frac{\1_{\{X_k=x-2\}}}{Z_k(x-2)^{(\g+\e)^{-1}}}\le\frac{\e}{2}
\end{align*}
again if $Z_m(x-2)\ge\Cst(C,\g,\e)$.

Third, 
\begin{equation}
\label{t3}
\Pb(T_3<T_1\wedge T_2|\F_m)=\Es[\Pb(T_3<T_1\wedge T_2|\F_S)\1_{\{S<\iy\}}],
\end{equation}
where 
$$S:=\inf\{n\ge m\tx{ s.t. }Z_n(x-2)\ge C Z_n(x)^{\g}\}.$$
By Proposition \ref{pol} {\bf(d)}, 
\begin{align}
\label{t4}
\Pb(T_3<T_1\wedge T_2|\F_S)\le\Cst(\e)[Z_S(x-2)^{-1}+Z_S(x)^{-1}]\le\Cst(C,\e)Z_m(x)^{-\g}.
\end{align}
Now inequalities \eqref{t1}, \eqref{t2}, \eqref{t3} and \eqref{t4} enable us to conclude.
\ep
\subsubsection{Proof of Proposition \ref{pemvolk}}
{\bf(a)}For all $n_0\ge\Cst$, $\e>0$, $\al\in(0,1)$, 
\bal
\Pb(\{Z_{n_0+x}(x-2)=Z_{n_0+x}(x+2)=2\}&\cap\{Z_{n_0+x}(x-1)\in((\al-\e)n_0/2,(\al+\e)n_0/2\}\\
&\cap\{Z_{n_0+x}(x)\in(1/2-\e,1/2+\e)n_0\})\ge\Cst(n_0).
\end{align*}
Now apply Corollary \ref{cor-loc}  twice, with $\pm:=+$ and $\pm:=-$, and Proposition \ref{bienvenue}. 

{\bf(b)} For all $x\le X_0$, let us prove that $\Pb(x-2\not\in R|x\in R)\ge\e>0$, which will imply the conclusion. For all $n_0\in\N$, let $u_{n_0}(x)$ be the first time $t$ such that $Z_t(x)=n_0$, and let $t(x)$ be the time of first visit to $x$. Then
\bal
\Pb(\{Z_{u_{n_0}(x)}(x-2)=Z_{u_{n_0}(x)}(x-1)=1\}&\cap\{Z_{u_{n_0}(x)}(x+1)\ge n_0\})\ge\Cst(n_0).
\end{align*}
Now apply Corollary \ref{cor-loc} with $x:=x$ and $\pm:=-$.

  \section{Rubin continuous time-lines construction}
 \label{rubin}
 \subsection{$W$-urn}
Again consider $W$-urn process with two colours $-1$ and $1$ studied in Section \ref{3v}, defined as follows: start with $Z_0(1)$ and  $Z_0(-1)$ balls of colours $1$ and $-1$ respectively. At each time step $n\ge0$, pick a ball of colour $i\in\{-1,1\}$ in the urn with a probability $W(Z_n(i))/(W(Z_n(-1))+W(Z_n(1)))$, and put it back  together with a ball of the same colour.

%%%%%%%%%%%%%%%%%%%%%%%FIGURE 1 %%%%%%%%%%%%%%%%%%%%%%%
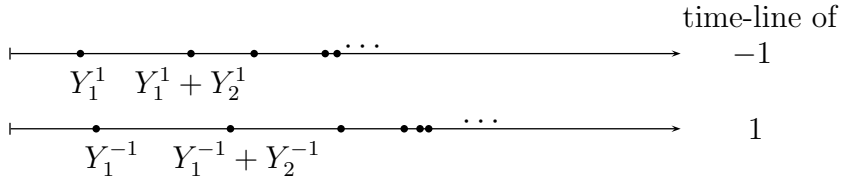
\psset{xunit= 1.05cm,yunit= 1cm}
%\psset{xunit= 0.4725cm,yunit= 0.45cm}
\pspicture(0,-1)(15,2)
\psset{linewidth=.5pt}
%draw time lines
\psline{|->}(0,0)(8.5,0)
\psline{|->}(0,1)(8.5,1)

\pscircle*(1.1,0) {.05} 
\pscircle*(2.8,0) {.05} 
\pscircle*(4.2,0){.05}
\pscircle*(5,0){.05}
\pscircle*(5.2,0){.05}
\pscircle*(5.31,0){.05}

\pscircle*(0.9,1){.05} 
\pscircle*(2.3,1){.05} 
\pscircle*(3.1,1){.05}
\pscircle*(4,1){.05}
\pscircle*(4.15,1){.05}
%text
\rput[u](1,0.6){$Y_1^1$}
\rput[u](2.3,0.6){$Y_1^{1}+Y_2^{1}$}
\rput[u](1.3,-.4){$Y_1^{-1}$}
\rput[u](3,-.4){$Y_1^{-1}+Y_2^{-1}$}
%ldots
\rput[b](4.5,1){$\cdots$}
\rput[b](6,0){$\cdots$}
%more text
\rput[u](9.5,1.5){time-line of}
\rput[l](9.15,1){$-1$}
\rput[l](9.35,0){$1$}
\endpspicture
\medskip

We now construct a continuous-time process $(\tZ_t(1),\tZ_t(-1))_{t\in\R_+}$ taking values in $\N^2$, which will be equal in law to $(Z_n(1),Z_n(-1))_{n\in\N}$, seen from the times of jumps. To this end we add balls of colour $i\in\{-1,1\}$ at rate $W(Z_n(i))$, using the following time-lines construction:

\begin{itemize}
\item
Let $(Y_k^{-1})_{k\in\N}$ and $(Y_k^{1})_{k\in\N}$ be two sequences of independent random variables of exponential law with $\Es Y_k^{i}=W(Z_0(i)+k-1)^{-1}$, independent from each other.

\item
Each of the colours $1$ and $-1$ has a clock with an alarm, set initially to $Y_1^{1}$ and $Y_1^{-1}$ respectively.

\item
Each time an alarm rings, we add a ball of the corresponding colour, say $i$. The other clock $-i$ keeps running, while the new alarm with $i$ is set at at time distance $Y_{k+1}^{i}$ if $k$ balls of colour $i$ have already been added in the urn.
\end{itemize}

More precisely let, for all $i\in\{-1,1\}$ and $t\in\R_+$, let
$$S_1:=\left\{\sum_{k=1}^n Y_k^{i},\,n\in\N\right\},\,S:=S_1\cup S_{-1},$$
let $\xi_n$ be the $n$-th smallest element in $S$, with the convention that $\xi_0:=0$.
and let
\bal
\tZ_t(i)&:=\sup\left\{l\ge0\tx{ s.t. }\sum_{k=1}^lY_k^{i}\le t\right\}+Z_0(i).
\end{align*}
\bl
\label{wscd}
The processes $(\tZ_{\xi_n}(1),\tZ_{\xi_n}(-1))_{n\ge0}$ and $(Z_n(1),Z_n(-1))_{n\ge0}$ (from a $W$-urn) are equal in law. 
\el
\bp
By memoryless property of exponentials, for all $n\in\N$, the joint law of the first alarms after time $\xi_n$ in the time-lines $-1$ and $1$ (conditioned on their past up to that time) is a couple of independent random variables of parameters $W(\tZ_{\xi_n}(1))$ and  $W(\tZ_{\xi_n}(-1))$ respectively. 

Now, if $U$ and $V$ are two independent random variables of parameters $u$ and $v$, then $\Pb[U<V]=u/(u+v)$, which completes the proof.
\ep

This result enables us to conclude that, if $W$ is reciprocally summable, then only one of the balls is taken in the urn infinitely often.
\begin{prop}
\label{locurn}
If $\sum_{k\in\N}1/W(k)<\iy$, then $Z_\iy(1)<\iy$ or $Z_\iy(-1)<\iy$ a.s. 
\end{prop}  
\bp
We have
\bal
\Pb\left(Z_\iy(1)=Z_\iy(-1)=\iy\right)&=\Pb\left(\sum_{k=1}^\iy Y_k^1=\sum_{k=1}^\iy Y_k^{-1}\right)\\
&=\Pb\left(Y_1^1=\sum_{k=1}^\iy Y_k^{-1}-\sum_{k=2}^\iy Y_k^1\right)=0,
\end{align*}
since $Y_1^1$ has continuous density, independent from $\sum_{k=1}^\iy Y_k^{-1}-\sum_{k=2}^\iy Y_k^1$.
\ep
\subsection{ESIRW on a locally finite graph $(G,\sim)$}
\label{erubin}
Let us similarly construct a continuous-time process $(\tX_t)_{t\in\R_+}$, which will be equal in law to $(X_n)_{n\ge0}$ for ESIRW, seen from times of jumps. Let $E(G)$ be the set of (non-oriented) edges of $(G,\sim)$. The process starts at $X_0:=x_0$ at time $0$: 

%%%%%%%%%%%%%%%%%%FIGURE 2%%%%%%%%%%%%%%%%%%%%%%%%%
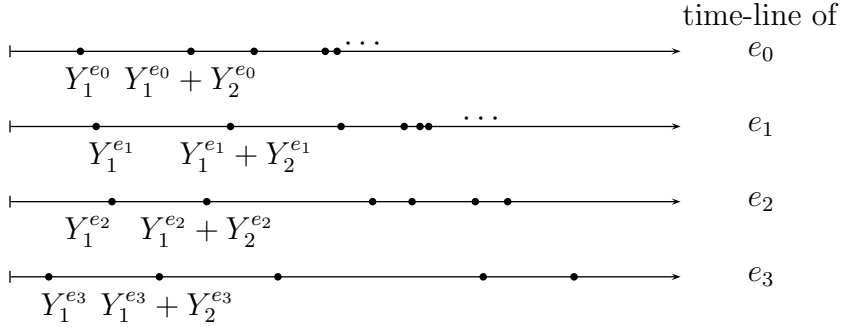
\psset{xunit= 1.05cm,yunit= 1cm}
%\psset{xunit= 0.4725cm,yunit= 0.45cm}
\pspicture(0,-2.5)(15,2)
\psset{linewidth=.5pt}
%draw time lines
\psline{|->}(0,0)(8.5,0)
\psline{|->}(0,1)(8.5,1)
\psline{|->}(0,-1)(8.5,-1)
\psline{|->}(0,-2)(8.5,-2)

\pscircle*(1.1,0) {.05} 
\pscircle*(2.8,0) {.05} 
\pscircle*(4.2,0){.05}
\pscircle*(5,0){.05}
\pscircle*(5.2,0){.05}
\pscircle*(5.31,0){.05}

\pscircle*(0.9,1){.05} 
\pscircle*(2.3,1){.05} 
\pscircle*(3.1,1){.05}
\pscircle*(4,1){.05}
\pscircle*(4.15,1){.05}

\pscircle*(1.3,-1) {.05} 
\pscircle*(2.5,-1) {.05} 
\pscircle*(4.6,-1){.05}
\pscircle*(5.1,-1){.05}
\pscircle*(5.9,-1){.05}
\pscircle*(6.31,-1){.05}

\pscircle*(0.5,-2){.05} 
\pscircle*(1.9,-2){.05} 
\pscircle*(3.4,-2){.05}
\pscircle*(6,-2){.05}
\pscircle*(7.15,-2){.05}

%text
\rput[u](1,0.6){$Y_1^{e_0}$}
\rput[u](2.3,0.6){$Y_1^{e_0}+Y_2^{e_0}$}
\rput[u](1.3,-.4){$Y_1^{e_1}$}
\rput[u](3,-.4){$Y_1^{e_1}+Y_2^{e_1}$}
\rput[u](1,-1.4){$Y_1^{e_2}$}
\rput[u](2.5,-1.4){$Y_1^{e_2}+Y_2^{e_2}$}
\rput[u](0.7,-2.4){$Y_1^{e_3}$}
\rput[u](2,-2.4){$Y_1^{e_3}+Y_2^{e_3}$}
%ldots
\rput[b](4.5,1){$\cdots$}
\rput[b](6,0){$\cdots$}
%more text
\rput[u](9.5,1.5){time-line of}
\rput[l](9.35,1){$e_0$}
\rput[l](9.35,0){$e_1$}
\rput[l](9.35,-1){$e_2$}
\rput[l](9.35,-2){$e_3$}
\endpspicture
\medskip

\begin{itemize}
\item
Let $(Y_k^{e})_{e\in E(G),k\in\N}$ be a collection of independent exponential random variables with $\Es Y_k^{e}=W(k-1)^{-1}$.

\item
Each edge $e$ has its own clock, which only runs when the process $(\tX_t)_{t\ge0}$ is adjacent to $e$.

\item Each time an edge $e$ has just been crossed, and at time $0$, its clock sets up an alarm at distance $Y_{k+1}^{e}$ if $e$ has been crossed $k$ times so far ($Y_1^{e}$ at time $0$).  

\item
Each time an edge $e$ sounds an alarm, $\tX_t$ crosses it instantaneously.  
\end{itemize}

Let $\tau_n$ be the $n$-th jump time of $(\tX_t)_{t\ge0}$, with the convention that $\tau_0:=0$.
\bl
\label{icd} (Davis \cite{davis}, Sellke \cite{sellke2})
The processes $(\tX_{\tau_n})_{n\ge0}$ and $(X_n)_{n\ge0}$ have the same distribution.
\el
The proof of Lemma \ref{icd} is left to the reader, being similar to that of  \ref{wscd}.

Let
$$\Gg_\iy:=\{e\in E(G)\tx{ s.t. }Z_\iy(e)=\iy\}.$$
\begin{prop}
If $\sum_{n\in\N}1/W(n)<\iy$, then $\Gg_\iy$ contains no even cycle.
\end{prop}
\bp
For simplicity, let us denote an even cycle by $\Z/l\Z$, $l$ even. Let, for all $i\in\Z/l\Z$, 
$$T^i:=\sum_{k=1}^\iy Y_k^{\{i,i+1\}}.$$
Then 
$$\left\{\Z/l\Z\subset\Gg_\iy\right\}\subset\left\{\sum_{x\in\Z/l\Z}(-1)^xT^x=0\right\}.$$
Now $\sum_{x\in\Z/l\Z}(-1)^xT^x\ne0$ a.s., which implies that 
$$\Pb\left(\Z/l\Z\subset\Gg_\iy\right)=0.$$
\ep

The technique carries over to show \cite{sellke2,limic} that, on graphs on bounded degree and if $W$ is reciprocally summable, then $\Gg_\iy$ is either a single edge or an odd cycle. 
\iffalse 
\subsection{VSIRW on a locally finite graph $(G,\sim)$}
\label{vrubin}
The time-lines construction can be easily be adapted \cite{bienvenue} to that case, by putting the alarms on vertices instead of edges, through a collection $(T_i^v)_{v\in G, i\in\N}$ of exponential independent random variables with expectation $\Es T_k^v=W(k)^{-1}$. 

\iffalse
It is natural to expect that, if $W$ is reciprocally summable, the walk would eventually localize on a single edge, at least with a regularity assumption on $W$, as for ESIRWs  \cite{limic-tarres}, but the results on that topic so far are only able to reject certain subsets as possible 
\fi
As for ESIRWs, the construction 
\fi
\section{Short proof of a.s. localization of the VRRW on $\Z$ on five consecutive sites}
 \label{short}
  Let, for all $x\in\Z$, 
 $$\Om(x)=\{\inf R'=x\}. $$
 Then
 $$\Pb\left(\cup_{x\in\Z}\Om(x)\right)=1,$$ 
 by Proposition \ref{pemvolk} {\bf(b)}, which asserts that the walk a.s. localizes on finitely many vertices. 
 The aim of this section is to prove the following two propositions. 
 \begin{prop} 
 \label{propexc}
 For all $x\in\Z$, $\Ep(x)\subset\{Z_\iy(x-1)<\iy\}\cup\{Z_\iy(x+1)<\iy\}$ a.s.
 \end{prop}
 \begin{prop}
 \label{omep}
For all $x\in\Z$, $\Om(x)\subset\Ep(x+4)$ a.s.
 \end{prop}
 These will imply Theorem \ref{loc5}, i.e. a.s. localization on the VRRW on five consecutive vertices: a.s. on $\Om(x)$, $Z_\iy(x+3)<\iy$ or $Z_\iy(x+5)<\iy$ by Propositions \ref{propexc} and \ref{omep}, and the former would imply that $R'=\{x,x+1,x+2\}$, which holds with probability $0$ by Proposition \ref{bienvenue} {\bf(b)}. 
 
 We first propose an alternative time-lines construction for VRRWs in Section \ref{arubin}, which will enable us to couple two random walks in Section \ref{coupling} with the following property: we will say that $\M'$ is greater than $\M$ if, at the time of $n$-th visit to any site $x\in\Z$, $\M'$ has more visited the right-hand side neighbour $x+1$ than $\M$, whereas on the contrary $\M$ has more visited $x-1$ than $\M'$. Then we will prove Propositions \ref{propexc} and \ref{omep} in Sections \ref{pexc} and \ref{pomep}.
 \subsection{An alternative time-lines construction for VSIRWs on $\Z$}
 \label{arubin}
We introduce the following time-lines construction on directed edges $\vec{E}(\Z)$ of $\Z$, which will enable us to introduce naturally a coupling in Section \ref{coupling}. 
\iffalse
Define 
 $$\vec{E}(\Z):=\{(x,y)\in G^2: x\sim y\}$$
 be the set of oriented edges naturally associated to the (nonoriented) tree $(G,\sim)$. 
 \fi
 If $e=(x,y)\in\vec{E}(\Z)$, let $\ue:=x$,
 $\ove:=y$, $\s(e):=(y,x)$. 
 
 The  continuous time process $(\tX_t)_{t\in\R_+}$ taking values in $\Z$ will be defined as follows: 
\begin{itemize}
\item
Let $(Y_k^{e})_{e\in \vec{E}(\Z),k\in\N}$ be a collection of independent exponential random variables with expectation one.

\item
Each oriented edge $e\in\vec{E}(\Z)$ has its own clock, which only runs when the process $(\tX_t)_{t\ge0}$ is adjacent to $e$.

\item Each time an edge $e$ has just been crossed, the clock of $\s(e)$ sets up an alarm at distance $Y_{k+1}^{\s(e)}/W(\tZ_t(\ue))$, if $\s(e)$ has been crossed $k$ times so far. At time $0$, we set up an initial alarm, at time distance $Y_1^{e}$,  for the edges $(x,x+1)$, $x\ge x_0$, $(x,x-1)$, $x\le x_0$.
\item
Each time an edge $e$ sounds an alarm, $\tX_t$ crosses it instantaneously.  
\end{itemize}
Let $\tau_n$ be the $n$-th jump time of $(\tX_t)_{t\ge0}$, with the convention that $\tau_0:=0$.
\bl
\label{aicd} 
The processes $(\tX_{\tau_n})_{n\ge0}$ and $(X_n)_{n\ge0}$ have the same distribution.
\el
The proof of Lemma \ref{aicd} is again left to the reader.
 \subsection{Coupling}
 \label{coupling}
Let us denote by $\M$ the function which maps a (deterministic) ``collection of alarms'' 
 $\Y=(Y_k^{e})_{e\in \vec{E}(\Z),k\in\N}$ and an initial site $x_0$ to a continuous-time (deterministic) walk $\M(\Y,x_0)$ on the vertices of $\Z$, as prescribed in Section \ref{arubin}.
 \begin{defin}
 Given $\Y=(Y_k^{e})_{e\in \vec{E}(\Z),k\in\N}$ and $\Y'=((Y')_k^{e})_{e\in \vec{E}(\Z),k\in\N}$ two collections of random variables on $\R_+$, we say that
 $\Y'\gg\Y$ if, for all $k\in\N$, $x\in\Z$, $(Y')_k^{(x,x+1)}\le Y_k^{(x,x+1)}$ and 
 $(Y')_k^{(x,x-1)}\ge Y_k^{(x,x-1)}$ a.s.
 \end{defin}
 
 Given $\Y=(Y_k^{e})_{e\in \vec{E}(\Z),k\in\N}$ and $\Y'=((Y')_k^{e})_{e\in \vec{E}(\Z),k\in\N}$ two collections of random variables on $\R_+$ we let, by a slight abuse of notation, $\M=(\tX_t)_{t\in\R_+}:=\M(\Y,x_0)$  and  $\M'=(\tX'_t)_{t\in\R_+}:=\M(\Y',x_0)$ be the continuous-time random walks starting associated to $\Y'$ and $\Y'$.
  
For all $i\in\N$, $u>0$, $j\in\Z$ and $e\in E(\Z)$ (resp. $\vec{e}\in\vec{E}(\Z)$) non-oriented (resp. oriented) edge, let $n_{e}(i)$ be the $i$-th visit time to $e$, let $l_j(t)$ be the local time at $j$ at time $t$, let $t_j(u):=\inf\{t\ge 0\tx{ s.t. }l_j(t)=u\}$, and let $T_j$ be the total time spent in $j$ for the random walk $\M$; let $n'_{e}(i)$ and $T'_j$ be the similar notation for $\M'$.

For any non-oriented edge $e=\{j,j+1\}$, let $\ove:=j+1$ and $\ue:=j$.

\begin{defin}
For all $i\in\N$ and $e\in E(\Z)$, let us define the property 
$E_{i,e}$ as follows:
$$Z'_{n'_{e}(i)}(\ove)\ge Z_{n_{e}(i)}(\ove) 
\text{ and }Z'_{n'_{e}(i)}(\ue)\le Z_{n_{e}(i)}(\ue),$$ 
with the convention that  $E_{i,j}$ holds whenever
$n_{e}(i)=\infty$ or $n'_{e}(i)=\infty$.
\end{defin}
\bl
\label{pcoupling}
Assume $\Y'\gg\Y$, and $W$ is nondecreasing. Then, for all $i\in\N$ and $e\in E(\Z)$, $E_{i,e}$ holds a.s. 
\el
\bp
Let, for all $T>0$,
$$\Pc_T:=\{e\in E(\Z), \, i\in\N\tx{ s.t. } n_e(i)\le T\tx{ and }n'_e(i)\le T, E_{i,e}\tx{ holds}\}.$$
Note that the property $\Pc_T$ can only change on a discrete set of times; we prove it by induction. Assume that $\Pc_T^-$ holds, i.e. that $\Pc_t$ holds for all $t<T$. We want to deduce $\Pc_T$: assume for instance that $n_e(i)=T$, $n'_e(i)\le T$, with $e=\{j,j+1\}$, $j\ge x_0$, and $i$ odd, so that $\tX_{n_e(i)}=\tX'_{n'_e(i)}=j+1$ (the other cases are similar). Obviously, 
$$Z'_{n'_e(i)}(j+1)=Z'_{n'_e(i-1)}(j+1)+1\ge Z_{n_e(i)}(j+1)=Z_{n_e(i-1)}(j+1)+1.$$
It remains to prove that 
\begin{equation}
Z'_{n'_{e}(i)}(j)\le Z_{n_{e}(i)}(j).
\end{equation}
Let $l_j:=l_j(n_i(e))-l_j(n_{i-1}(e))$ (resp. $l'_j$) be the local time spent at $j$ between times $n_e(i-1)$ and $n_e(i)$ (resp. $n'_e(i-1)$ and $n'_e(i)$). Then
$$l_j=\frac{Y_{(i+1)/2}^{(j,j+1)}}{W(Z_{n_e(i-1)}(j+1))}\ge\frac{(Y')_{(i+1)/2}^{(j,j+1)}}{W(Z'_{n'_e(i-1)}(j+1))}=l'_j.$$
Let
$$u:=\inf\{0\le s\le l'_j\tx{ s.t. } Z'_{t'_j(n'_e(i)+s)}(j)> Z_{t_j(n_e(i)+s)}(j)\}.$$
Let us now consider the last jump occuring strictly before the local time  at site $j$ is $u$: it has to be a move from $\M'$, from $j$ and $j-1$ and back (possibly with a simultaneous move from $\M$). At that last time, the numbers of visits to $(j,j-1)$ are equal for $\M$ and $\M'$ (since those of $j$ and $(j,j+1)$ are) and, by $\Pc_T^-$, the numbers of visits to $j-1$ is greater for $\M$ than for $\M'$. Consequently $\M$ must move before $\M'$, in $j$-th local time, after that last jump before $u$. Therefore $u>l'_j$, and $\Pc_T$ holds. 
\ep
 \subsection{Proof of Proposition \ref{propexc}}
 \label{pexc}
Fix $x\in\Z$. 
Let $\Y:=(Y_k^{e})_{e\in \vec{E}(\Z),k\in\N}$ be a collection of independent exponential random variables with expectation $1$, and let 
$$\Y^{'(n)}:=((Y^{'(n)})_k^{e})_{e\in \vec{E}(\Z),k\in\N}:=(Y_k^{e}+\1_{\{e=(x,x-1)\}}\1_{\{k=n\}})_{e\in \vec{E}(\Z),k\in\N}.$$
Let  $\M=(\tX_t)_{t\in\R_+}:=\M(\Y,x_0)$  and  $\M^{'(n)}=(\tX'_t)_{t\in\R_+}:=\M(\Y',x_0)$, and let us use the notation from Section \ref{coupling}.

Let
\bal
\Qc&:=\{Z_\iy(x+1)=Z_\iy(x-1)=\iy\}\cap\{T_x<\iy\},\\
\Qc^{'(n)}&:=\{Z'_\iy(x+1)=Z'_\iy(x-1)=\iy\}\cap\{T'_x<\iy\}.\\
\end{align*}
\bl
\label{exc}
For all $n\in\N$, $\Pb(\Qc\cap\Qc^{'(n)})=0$.
\el
\bp
If $Z_\iy(x+1)=Z_\iy(x-1)=Z'_\iy(x+1)=Z'_\iy(x-1)=\iy$ and $T_x<\iy$, $T'_x<\iy$, then
$$T_x=\sum_{k=1}^n\frac{Y_k^{(x,x+1)}}{Z_{n_{(x+1,x)}(k)}(x+1)}=
\sum_{k=1}^n\frac{Y_k^{(x,x-1)}}{Z_{n_{(x-1,x)}(k)}(x-1)}$$
and, using Lemma \ref{pcoupling},
$$T'_x=\sum_{k=1}^n\frac{(Y')_k^{(x,x+1)}}{Z'_{n'_{(x+1,x)}(k)}(x+1)}\le T_x<
\sum_{k=1}^n\frac{(Y')_k^{(x,x-1)}}{Z'_{n'_{(x-1,x)}(k)}(x-1)}=T'_x\tx{ a.s.,}$$
which is contradictory.
\ep

Let $\F_n:=\s(\tX_0,\ldots,\tX_n)\subset\s(Y_k^{e},1\le k\le n)$. Then Lemma \ref{exc} implies
$$\Pb(\Qc^c\cup(\Qc^{'(n+1)})^c|\F_n)=1.$$
But
$$\Pb(\Qc^c|\F_n)\ge e^{-1}\Pb((\Qc^{'(n+1)})^c|\F_n),$$
so that 
$$\Pb(\Qc^c|\F_n)\ge (1+e)^{-1}.$$
Now $\Pb(\Qc^c|\F_n)\des_{n\to\iy}\1_{\Qc^c}$ a.s., so that $\Qc^c$ holds almost surely.
 \subsection{Proof of Proposition \ref{omep}} 
 \label{pomep}
 By Corollary \ref{noinf} {\bf(a)-(b)}, for all $x\in\Z$,
 $$\Om(x)\subset\Ep(x)\subset\Ep(x)\cap\{\al_n^-(x+2)\des_{n\to\iy}\al_\iy^-(x+2)\in[0,1)\}
 \subset\Ep(x+4)\cup\{\al_\iy^-(x+2)=0\}. $$
 Therefore, if we let 
 $$\ov{\Om}(x):=\Om(x)\cap\{\al_\iy^-(x+2)=0\},$$
 it is sufficient to show that $\Pb(\ov{\Om}(x))=0$.
 
 Let us assume $x:=0$ for simplicity. Our proof relies on the study of the asymptotic behaviour of $\log Z_n(3)/Z_n(2)$: roughly speaking, we will show that this quantity must converge to $0$ on $\ov{\Om}(0)$ but that, on the other hand, its convergence to $0$ can only happen with zero probability, due to unstability.
 
 For all $n\in\N$ and  $x\in\Z$, let
 \bal
R_n&:=Z_n(4)+Z_n(2)-Z_n(1)-Z_n(3),\\
z_n&:=\log\frac{Z_n(3)}{Z_n(2)},\,\,\,y_n:=\frac{R_n}{Z_n(2)Z_n(3)},
 \end{align*}
and let $t_n(x)$ (resp. $t_n^\pm(x)$) be the $n$-th visit time of $x$ (resp. $(x,x\pm1)$, counted once the edge has been visited), and let $Z_n^\pm(x)$ be the number of visits of the edge $(x,x\pm1)$ at time $n$.

Note that, for all $i\in\N$, 
\begin{equation}
\label{ri}
R_i=Z_i^-(5)-Z_i^+(0)+\1_{\{X_i=2\tx{ or }X_i\ge4\}}+\Cst(x_0),
\end{equation}
so that
\begin{equation}
\label{ri2}
R_{t_i^+(2)}=Z_{t_i^+(2)}^+(4)-Z_{t_i^+(2)}^-(1)+\Cst(x_0).
\end{equation}
For all $\al\in(0,1)$ and $k\ge j$, we write $j\lr_\al k$ if, for all $i\in[j,k]$, $\al_i^-(2)\vee\al_i^+(3)\le\al$. 
 \bl
 \label{evq}
 Assume $Z_n(2)\ge\Cst$. Then there exist $\al_0:=\Cst$, $(\eta_{j,k})_{k\ge j\ge n}$ and $(r_{j,k})_{k\ge j\ge n}$ such that,    and $k\ge j\ge n$ with $j\lr_{\al_0} k$,   
 \begin{equation}
 \label{eevq}
 z_k-z_j
 =\sum_{i=j+1}^k {\textnormal{\1}}_{\{X_{i-1}=2,X_i=3\}}y_{i}+\eta_{n,k}+r_{n,k},
 \end{equation}
 where, for all $\al\le\al_0$ and $\e>0$,
\begin{equation}
\label{estevq}
 \Pb\left(\sup_{k\ge j\ge n}|\eta_{j,k}|{\textnormal{\1}}_{j\lr_\al k}\ge\e|\F_n\right)\le\Cst\frac{\al}{\e^2Z_n(2)},\,\,\,\,
 |r_{j,k}|\le\frac{\Cst}{Z_j(2)}.
 \end{equation}
  \el
 Lemma \ref{evq} is proved in Section \ref{pevq}. 
 \bl
  \label{omest}
  $\ov{\Om}(0)\subset\{\limsup(\al_n^+(3)/\al_n^-(2))\le1\}\cap\{\sum|y_{t_n^+(2)}|{\textnormal{\1}}_{t_n^+(2)<\iy}<\iy\}$.
  \el
  Lemma \ref{omest} is proved in Section \ref{pomest}.
  
For simplicity we will use the notation $t_n:=t_n^+(2)$ until the beginning of Section \ref{pevq}.
Let $E_1:=\{(2,1),(0,1)\}$ be the set of edges pointing to $1$.  We introduce the continuous-time construction of Section \ref{arubin} in order to analyze the walk after time ${t_n}$, but modify the rule for the two edges in $E_1$: each time an edge $e\in\s(E_1)$ has been crossed, we set up an alarm on $\s(e)$ at time distance $V_{k+1}$ if site $1$ has been visited $k$ times after time $n$. The walk thus defined depends on the sequences $(Y_k^{e})_{e\in\vec{E}(\Z)\setminus E_1}$, $(V_k)_{k\in\N}$ and on the starting point $X_{t_n}$.  
  
Let  $(Y_k^{e})_{e\in\vec{E}(\Z)\setminus E_1}$ (resp. $((Y')_k^{e})_{e\in\vec{E}(\Z)\setminus E_1}$) be collections of independent exponential random variables with expectation one, and let $(V_k)_{k\in\N}$ (resp. $(V'_k)_{k\in\N}$) be sequences of independent exponential random variables with expectation $\Es\, V_k=(n_0+k-1)^{-1}$ (resp. $\Es\, V'_k=(n_0+k-1)^{-1}$); these define two continuous-time walks $\M$ and $\M'$. We choose $n_0:=Z_{t_n}(1)$ and 
$n'_0:=Z_{t_n}(1)-a\sqrt{Z_{t_n}(1)}$. 

Note that both $\M$ and $\M'$ are continuous-time versions of a VRRW: $\M$ of the original walk $(X_k)_{k\ge {t_n}}$, and $\M'$ of the VRRW $:=(X'_k)_{k\ge {t_n}}$ taken after discounting the number of visits to $1$ of $a\sqrt{Z_{t_n}(1)}$. 

Given $\e>0$ and $n\in\N$, assume $\al_{t_n}^-(2)<\e$, and let
\bal
 T_n:=\inf\{k\ge t_n  &\tx{ s.t. }\al_k^+(3)\vee\al_k^-(2)\ge(1+\e)\al_n^-(2)\\
 &\tx{ or }\sum_{i=k+1}^n|y_{i-1}|\1_{\{X_{i-1}=2,X_i=3\}}>\e\tx{ or }X_k=0\}.
 \end{align*}
 
We let $Z'_.(.)$, $z'$, $R'_.$, $T'$, $\ov{\Om}^{\,'}(0)$ and $t'_.$ be the notation for $\M'$ equivalent to the one already defined for $\M$.

Let 
\bal
\Qc_n&:=\ov{\Om}(0)\cap\{T_n=\iy\},\\
\Qc_n^{'}&:=\ov{\Om}^{\,'}(0)\cap\{T'_n=\iy\}.\\
\end{align*}
\bl
\label{0int}
If $\e\le\Cst$, $Z_{t_n}(1)\ge\Cst$ and $a\ge\Cst$, then $\Pb((\Qc_n\cap\Qc_n^{'})^c\,|\,\F_{t_n})\ge\Cst$.
\el
Lemma \ref{0int} is proved in Section \ref{s0int}.
\bl 
\label{prime}
For all $\de>0$, there exists $C(a,\de)$ (depending only on $a$ and $\de$) such that 
$$\Pb((\Qc'_n)^c\,|\,\F_{t_n})\le C(a,\de)\Pb(\Qc^c_n\,|\,\F_{t_n})+\de.$$
\el
Lemma \ref{prime} completes the proof of Proposition \ref{pomep}: indeed, using also Lemma \ref{0int}, for all $n\in\N$, 
$$\Cst\le\Pb(\Qc_n^c\,|\,\F_{t_n})+\Pb((\Qc_n^{'})^c\,|\,\F_{t_n})\le(1+C(a,\de))\Pb(\Qc^c_n\,|\,\F_{t_n})+\de,$$
which implies $\Pb(\Qc^c_n\,|\,\F_{t_n})\ge\Cst(a)$ if we choose $\de\le\Cst$. Now $\Pb(\Qc^c_n\,|\,\F_{t_n})\des_{n\to\iy}\1_{\Qc^c_n}$ a.s., so that $\Qc^c_n$ holds a.s. Now $\liminf_{n\to\iy}\Qc^c_n=\ov{\Om}(0)^c$, since $T_n=\iy$ for some $n$ on $\ov{\Om}(0)$, using Lemma \ref{omest} and Corollary \ref{noinf} {\bf(a)}.

\bp
There is a one-to-one correspondence between $(V_k)_{k\in\Z}$ and a simple birth process $\{N_t,t\ge0\}$ with initial population size $n_0$, defined by 
$$N_t:=n_0+\sup\left\{k\in\N\tx{ s.t. }\sum_{i=1}^k V_i\le t\right\}.$$
By a result of D. Kendall \cite{kendall}, $\{N_{\log(1+t/W)}, t\ge0\}$ is a Poisson process with unit parameter, where $W:=\lim N_te^{-t}$ is a Gamma random variable $\G(n_0,1)$ with shape $n_0$ and scale $1$. 

The same remark applies to $(V'_k)_{k\in\Z}$. Now recall that $U\sim\G(\la,1)$ has density $\phi_\la(x):=x^{\la-1}e^{-x}/\G(\la)$ w.r.t. Lebesgue measure $\Lc(dx)$, expectation and variance $\la$. Assume $\la\ge1$; for all $a$, $c$ $>0$, there exist $c_1$, $c_2$ $:=\Cst(a,c)$ such that, for all $x$ such that $|x-\la|\le c\sqrt{\la}$, $\phi_\la(x)/\phi_{\la-\sqrt{\la}a}(x)\in[c_1,c_2]$. Hence, for any Borel subset $A$ of $\R$, 
\bal
\Pb(W'\in A)&\le\Pb(|W'-n'_0|\ge\sqrt{n'_0}c)+\Pb(W'\in A\cap[n'_0-\sqrt{n'_0}c,n'_0+\sqrt{n'_0}c])\\
&\le c^{-2}+\Cst(a,b)\Pb(W\in A).
\end{align*}
We conclude by choosing $c:=\de^{-1/2}$.
\ep
\subsubsection{Proof of Lemma \ref{evq}}
  \label{pevq}
  Assume $j\lr_{\al_0} k$, with $\al_0\le\Cst$.
 Let us first estimate $h(Z_k(3))-h(Z_j(3))$: using Proposition \ref{pol} {\bf(a)}, there exists a constant $C>0$ such that
 \bal
 h(Z_k(3))=Y_k^+(2)+Y_k^-(4)=Y_k^+(2)+U_{k,+}^+ -\frac{\1_{\{X_k\ge4\}}}{Z_{k-1}(3)}+C
 \end{align*}
But
$$U_{k,+}^+(3)=U_{k,+}^-(3)+\cU_{k,+}(3).$$
Now we estimate $U_{k,+}^-(3)$. We assume $X_0\le 2$ (the other case is similar): then, for all $i\in\N$, $t_i^+(2)<t_i^-(3)<t_{i+1}^+(2)$, and $Z_{t_i^-(3)}(x)=Z_{t_{i+1}^+(2)}(x)$ for $x=3,4$, so that 

\bal
&U_{k,+}^-(3)-U_{j,+}^-(3)=\sum_{i=Z_n^-(3)+1}^{Z_k^-(3)}\frac{Z_{t_i^-(3)}(4)}{Z_{t_i^-(3)}(3)Z_{t_i^-(3)}(2)}\\&
=\sum_{i=Z_j^-(3)+1}^{Z_k^-(3)}\left(\frac{Z_{t_{i+1}^+(2)}(4)}{Z_{t_{i+1}^+(2)}(3)Z_{t_{i+1}^+(2)}(2)}
+\frac{Z_{t_{i+1}^+(2)}(4)}{Z_{t_{i+1}^+(2)}(3)}\left(\frac{1}{Z_{t_{i}^-(3)}(2)}-\frac{1}{Z_{t_{i+1}^+(2)}(2)}\right)\right)\\
&=\sum_{l=j+1}^k\frac{\1_{\{X_{l-1}=2,X_l=3\}}}{Z_{l-1}(2)}\frac{Z_{l-1}(4)}{Z_{l-1}(3)}+r_{j,k}^1,
\end{align*}
where
$$|r_{j,k}^1|\le2\sup_{j\le l<k}\left(\frac{Z_{l}(4)}{Z_{l}(3)}\right)\frac{1}{Z_j(3)}\le\frac{1}{2Z_j(3)}$$
if $\al\le\Cst$. 

In summary, 
\bal
 \log\frac{Z_k(3)}{Z_j(3)}=\sum_{l=j+1}^k\frac{\1_{\{X_{l-1}=2,X_l=3\}}}{Z_{l-1}(2)Z_{l-1}(3)}(Z_{l-1}(2)+Z_{l-1}(4))
 +\cU_{k,+}(3)-\cU_{j,+}(3)+r_{j,k}^2,
\end{align*}
where $|r_{j,k}^2|\le3/Z_j(2)$, using \eqref{logapprox}.

Similarly, 
\bal
 \log\frac{Z_k(2)}{Z_j(2)}=\sum_{l=j+1}^k\frac{\1_{\{X_{l-1}=2,X_l=3\}}}{Z_{l-1}(2)Z_{l-1}(3)}(Z_{l-1}(1)+Z_{l-1}(3))
 +\cU_{k,-}(2)-\cU_{j,-}(2)+r_{j,k}^3,
\end{align*}
where $|r_{j,k}^3|\le3/Z_j(2)$. 

Now, if $\al_0\le\Cst$, then
$$\sum_{l=j}^k|y_{l-1}-y_l|\1_{\{X_{l-1}=2,X_l=3\}}\le\sum_{l=j}^k\frac{\1_{\{X_{l-1}=2,X_l=3\}}}{Z_i(2)Z_i(3)}\le\frac{2}{Z_j(2)}.$$

This provides \eqref{eevq}, letting
 $$\eta_{j,k}:=\cU_{k,+}(3)-\cU_{j,+}(3)-(\cU_{k,-}(2)-\cU_{j,-}(2)).$$

Let us now estimate $\eta_{j,k}$. The following Lemma \ref{estun} will enable us to conclude, using Chebyshev's inequality. 

For all $n\in\N$, $x\in\Z$ and $\al>0$, let
 \bal
 \cU_{n,\pm}(x)&:=\hU_{n,\pm}^+(x)-\hU_{n,\pm}^-(x)\\
 \tU_{n,\pm}(x)&:=\sum_{k=1}^n(\cU_{k,\pm}(x)-\cU_{k-1,\pm}(x))\1_{\{\al_{k-1}^\pm(x)\le\al\}}\\
 \end{align*}
  \bl
  \label{estun}
  For all $x\in\Z$, $k\ge n$, $\al\le\Cst$, 
  $$\Es\left(\sup_{k\ge n} \,(\tU_{k,\pm}(x)-\tU_{n,\pm}(x))^2|\F_n\right)\le\Cst\al Z_n(x)^{-1}.$$
  \el
  \bp
  $(\tU_{n,\pm}(x))_{n\ge0}$ is a martingale, and 
  $$\Es((\tU_{k+1,\pm}(x)-\tU_{k,\pm}(x))^2|\F_k)\le\frac{\1_{\{X_k=x\}}}{Z_k(x)^2}\left(\al+\frac{\al^2}{1-\al}\right)
  \le2\al\frac{\1_{\{X_k=x\}}}{Z_k(x)^2},$$
  which enables us to conclude by Doob's inequality.
  \ep

\subsubsection{Proof of Lemma \ref{omest}}
\label{pomest}
Given $\e>0$ and $n\in\N$, let 
$$T_n:=\inf\{k\ge n\tx{ s.t. } \al_k^-(2)\ge(1+\e^2)\al_n^-(2)\tx{ or }Y_k^+(0)-Y_n^+(0)>\e^2\}.$$
Let us first prove that,  for all $\e>0$ and sufficiently large $n\in\N$, if 
\begin{equation}
\label{hyp}
\al_n^+(3)\ge(1+\e)^4\al_n^-(2)\tx{ with }\al_n^-(2)\le\e^2\le\Cst,
\end{equation} 
then
\begin{equation}
\label{phalf}
\Pb(\{T_n<\iy\}\cup\{\liminf\al_k^-(3)>0\}|\F_n)\ge1/2.
\end{equation}
This will imply the first part of the inclusion. Indeed let, for all $p\in\N$, 
\bal
\A&:=\{\limsup_{k\to\iy}(\al_k^+(3)/\al_k^-(2))\le1\},\\
\B_p&:=\{T_p<\iy\}\cup\{\liminf_{k\to\iy}\al_k^+(3)>0\}\cup\{\limsup_{k\to\iy}\al_k^-(2)>0\}.
\end{align*}
Then \eqref{phalf} implies, for all $p\in\N$, that for large $n\in\N$, $1/2\le\Pb(\A\cup\B_p|\F_n)$, which converges a.s. to $\1_{\A\cup\B_p}$ as $n\to\iy$, so that $\A\cup\B_p$ holds a.s.; subsequently $\A\cup\liminf\B_p$ holds a.s. 
Now, by Proposition \ref{pol} {\bf(c)} and Corollary \ref{noinf} {\bf(a)} and {\bf(d)},
$$\ov{\Om}(0)\cap\{\liminf\al_k^+(3)>0\}\subset\ov{\Om}(0)\cap\Ep(-1,0,1)=\emptyset,$$ 
which implies $\ov{\Om}(0)\subset\limsup\B_p^c\subset\A$  and  enables to conclude. 

Let us assume \eqref{hyp}, and $\al_n^+(3)\le\e^2$ w.l.o.g (possibly by choosing $n$ larger), and show \eqref{phalf}.
Fix $\al:=(1+\e)^4\al_n^-(2)$, and let $p$ be the largest $i\in[n+1,k]$ such that $\al_{i-1}^-(2)>\al$ or $\al_{i-1}^+(3)>\al$; then $\al_p^+(3)\ge(1+\e)^2\al_p^-(2)$. Note that, for all $i\ge p$, $Z_i(2)/Z_i(3)\ge\al_i^-(3)\ge1-\e^2$, since $Z_i(3)\le Z_i(2)+Z_i(4)$.

Let $t_j:=t_j^+(2)$ for simplicity. In order to make use of Lemma \ref{evq}, we need to estimate $\sum_{j=Z_p^+(2)+1}^{Z_k^+(2)}y_{t_j}$. Now,  for all $i\ge p$, using \eqref{ri},
\begin{equation}
\label{230}
R_i\ge R_p-(Z_i^+(0)-Z_p^+(0))-1.
\end{equation}
On one hand, 
\begin{align}
\nonumber
\sum_{j=Z_p^+(2)+1}^{Z_k^+(2)}\frac{1}{Z_{t_j}(2)Z_{t_j}(3)}&\ge\sum_{j=0}^{Z_k^+(2)-Z_p^+(2)-1}\frac{1}{(Z_{p}(2)+j)(Z_{p}(2)+j+1)}(1-\Cst.\al)\\
\label{231}
&\ge\left(\frac{1}{Z_p(2)}-\frac{1}{Z_k(2)}\right)(1-\Cst.\al).
\end{align}
On the other hand, assume $Z_n(2)\ge\Cst$, so that $Z_i(2)Z_i(3)\ge(1-\e^2)^2Z_i(2)(Z_i(2)+1)$; then
\begin{align}
\nonumber
\De_{p,k}:=\sum_{i=p}^{k-1}\frac{Z_i^+(0)-Z_p^+(0)}{Z_i(2)Z_i(3)}\1_{\{X_i=2\}}&\le(1-\e^2)^{-2}\sum_{i=p}^{k-1}(Z_i^+(0)-Z_p^+(0))\left(\frac{1}{Z_{i-1}(2)}-\frac{1}{Z_{i}(2)}\right)\\
\label{232}
&\le(1-\e^2)^{-2}\sum_{i=p}^{k-1}\frac{\1_{\{X_{i-1}=0,X_i=1\}}}{Z_{i-1}(2)}\le(1-\e^2)^{-4}\al\e^2
\end{align}

Now, if $Z_n(1)\ge\Cst(\e)$, using \eqref{estevq},
$$\Pb(\sup_{k\ge j\ge n}|\eta_{j,k}|\1_{j\lr_\al k}\,|\,\F_n)\ge1-\frac{\Cst}{\e^4\al_n^-(2)}\left(\frac{1}{Z_n(2)}+\frac{1}{Z_n(3)}\right)\ge\frac{1}{2},$$
where we use $\al_n^-(2)Z_n(3)=\al_n^+(2)Z_n(1)\ge(1-\e^2)Z_n(1)$.

Now assume that $E$ holds.
If $Z_n(1)\ge\Cst(\e)$, then Lemma \ref{evq} implies, together with \eqref{231} and \eqref{232}, 
\bal
\log\frac{Z_k(3)}{Z_k(2)}-\log\frac{Z_p(3)}{Z_p(2)}+\Cst\e^2\al\ge R_p\left(\frac{1}{Z_p(3)}-\frac{1}{Z_k(3)}\right)
\ge\frac{R_p}{Z_p(3)}-\frac{R_k+Z_k(1)}{Z_k(3)},
\end{align*}
where we use in the second inequality that $R_k\ge R_p-Z_k(1)$, by \eqref{230}.

Therefore
\bal
\frac{Z_k(4)}{Z_k(2)}\ge\log\frac{Z_k(3)}{Z_k(2)}+\frac{R_k+Z_k(1)}{Z_k(3)}-\Cst\e^2\al\ge\frac{Z_p(4)-Z_p(1)}{Z_p(2)}-\Cst\e^2\al
\ge\frac{\e\al}{2}
\end{align*}

Therefore $\liminf\al_k^+(3)>0$ on $E$, which completes the proof of the first part of the inclusion.

Now assume $\ov{\Om}(0)\cap\{\limsup(\al_k^+(3)/\al_k^-(2))\le\}$ holds. Then $Z_n(3)/Z_n(2)$ converges to $1$, and the estimates \eqref{230}, \eqref{231} and \eqref{232} ensure that $\sum y_{i-1}^-\1_{\{X_{i-1}=2,X_i=3\}}<\iy$. On the other hand Lemma \ref{evq} ensures, together with Lemma \ref{estun} and Doob's martingale Theorem that 
$\limsup_{k,n\to\iy: k\ge n}|\sum y_{i-1}\1_{\{X_{i-1}=2,X_i=3\}}|<\iy$, which completes the second part of the inclusion.

\subsubsection{Proof of Lemma \ref{0int}}
\label{s0int}
 Let $\al:=(1+\e)\al_{t_n}^-(2)$, and let $k\ge j\ge n$ be such that $t_k\le T$, $t'_k\le T'$.
 
 Lemma \ref{evq} still holds for $\M'$, with the modification that $Z_i(1)$ has to be replaced by $Z_i(1)-a\sqrt{Z_{t_n}(1)}$ in $y_i$: letting $\eta'_{.,.}$, $\lr'_\al$ and $r'_{.,.}$ be the corresponding notation,
 $$z'_{t'_k}-z_{t_k}=z'_{t'_j}-z_{t_j}+\sum_{i=j+1}^k\frac{a\sqrt{Z_{t_n}(1)}}{Z'_{t'_i}(2)Z'_{t'_i}(3)}+\sum_{i=j+1}^k(y'_{t'_i}-y_{t_i})+\eta'_{t'_j,t'_k}-\eta_{t_j,t_k}+r'_{t'_j,t'_k}-r_{t_j,t_k}.$$
 Let us first estimate $y'_{t'_i}-y_{t_i}$, for all $i\in[n,k]$: the coupling of Lemma \ref{pcoupling} holds, as long as $n_e(i)\le T$, $n'_e(i)\le T'$. Therefore $Z'_{t'_i}(3)\ge Z_{t_i}(3)$ (and thus $z'_{t'_k}\ge z_{t_k}$), thus $Z^{'+}_{t'_i}(3)\ge Z^+_{t_i}(3)$ (since the numbers of visits of edge $(2,3)$ are equal for both walks at those times), which implies subsequently $Z^{'+}_{t'_i}(4)\ge Z^+_{t_i}(4)$. Similarly $Z^{'-}_{t'_i}(1)\le Z^-_{t_i}(1)$ so that, using identity \ref{ri2}, $R'_{t'_i}\ge R_{t_i}$. 
 
 On the other hand, using that $|x|\le2|\ln(1+x)|$ for $|x|<1$,
 $$\left|\frac{Z_{t_i}(2)Z_{t_i}(3)}{Z'_{t'_i}(2)Z'_{t'_i}(3)}
-1\right|
\le2\left|\ln\frac{Z_{t_i}(2)Z_{t_i}(3)}{Z'_{t'_i}(2)Z'_{t'_i}(3)}\right|
\le2\left(\ln\frac{Z'_{t'_i}(3)}{Z_{t_i}(3)}+\ln\frac{Z_{t_i}(2)}{
Z'_{t'_i}(2)}\right)=2(z'_{t'_i}-z_{t_i}).$$
Hence
\bal
y'_{t'_i}-y_{t_i}\ge-|R_{t_i}|\left|\frac{1}{Z_{t_i}(2)Z_{t_i}(3)}-\frac{1}{Z'_{t'_i}(2)Z'_{t'_i}(3)}\right|\ge-2|y_{t_i}|\,(z'_{t'_i}-z_{t_i}).
\end{align*}
Now let 
$$u(k):=\mathsf{argmax}_{i\in[n,k]}(z'_{t'_i}-z_{t_i}),~~\tau_k:=z'_{t'_{u(k)}}-z_{t_{u(k)}}.$$
Using \eqref{231} in the proof of Lemma \ref{omest}, and $\sum_{i=j+1}^k |y_{t_i}|\le\e$, we obtain that
\begin{align}
\nonumber
z'_{t'_k}-z_{t_k}\ge & \,\,\,z'_{t'_j}-z_{t_j}-2\e\tau_k+a\sqrt{Z_{t_n}(1)}(Z_{t_j}(2)^{-1}-Z_{t_k}(2)^{-1})(1-\Cst\al)\\
\label{variance}
&\,\,\,\,\,+\eta'_{t'_j,t'_k}-\eta_{t_j,t_k}+r'_{t'_j,t'_k}-r_{t_j,t_k}.
\end{align}
Given $b>0$, let
$$\De:=\left\{\sup_{k\ge j\ge n}|\eta_{t_j,t_k}|{\textnormal{\1}}_{j\lr_\al k}\le b\,\,\frac{\sqrt{Z_{t_n}(1)}}{Z_{t_n}(2)}\right\}\bigcap
\left\{\sup_{k\ge j\ge n}|\eta'_{t'_j,t'_k}|{\textnormal{\1}}_{j\lr'_\al k}\le b\,\,\frac{\sqrt{Z_{t_n}(1)}}{Z_{t_n}(2)}\right\}.$$
Then, by \eqref{estevq} and Chebyshev's inequality, if $b\ge\Cst$,  then
$$\Pb(\De\,|\,\F_{t_n})\ge1/2.$$ 
Now assume $\De$ holds, and apply \eqref{variance} for $j:=n$ and $k$ such that $Z_{t_k}(2)\ge 2Z_{t_n}(2)$: if  $Z_{t_n}(1)\ge\Cst(b)$, then
\begin{equation}
\label{intvar}
\tau_k\ge\left(\frac{a}{2}-3b\right)\frac{\sqrt{Z_{t_n}(1)}}{Z_{t_n}(2)}-2\e\tau_k.
\end{equation}
Apply again \eqref{variance} for $j:=u(k)$ and $k$, and use \eqref{intvar}:
$$z'_{t'_k}-z_{t_k}\ge (1-2\e)\tau_k-3b\frac{\sqrt{Z_{t_n}(1)}}{Z_{t_n}(2)}
\ge\left\{\left(\frac{1-2\e}{1+2\e}\right)\left(\frac{a}{2}-3b\right)-3b\right\}\frac{\sqrt{Z_{t_n}(1)}}{Z_{t_n}(2)}
\ge\frac{a}{3}\frac{\sqrt{Z_{t_n}(1)}}{Z_{t_n}(2)},$$
if $a\ge\Cst(b)$ and $\e\le\Cst$. 

We conclude by the remark that $z_{t_n}$ (resp. $z'_{t'_n}$) converge to $0$ on $\ov{\Om}(0)$ (resp. $\ov{\Om}^{\,'}(0)$), so that $\Pb(\Qc_n\cap\Qc'_n\cap\De)=0$, and $\Pb((\Qc_n\cap\Qc'_n)^c\,|\,\F_{t_n})\ge\Pb(\De\,|\,\F_{t_n})\ge1/2$ if $b\ge\Cst$.

 \section{Appendix}
 We recall Lévy's conditional Borel-Cantelli Lemma. Let $\Ff=(\F_n)_{n\in\N}$ be a filtration, and let $(\xi_n)_{n\in\N}$ be an $\Ff$-adapted sequence taking values in $\R_+$.
 \bl
 \label{bc}
 Assume $(\xi_n)_{n\in\N}$ is a.s. bounded by a constant $C>0$. Then
 $$\left\{\sum_{k=1}^\iy\xi_k<\iy\right\}=\left\{\sum_{k=1}^\iy\Es(\xi_k|\F_{k-1})<\iy\right\}.$$
 \el
 \bp
 Let 
 $$M_n:=\sum_{k=1}^n(\xi_k-\Es(\xi_k|\F_{k-1})),$$
with $\F_0:=\emptyset$. Then $(M_n)_{n\in\N}$ is a martingale, and 
 $$<M>_n:=\sum_{k=1}^{n-1}\Es((M_k-M_{k-1})^2|\F_{k-1})\le\sum_{k=1}^{n-1}\Es(\xi_k^2|\F_{k-1})=O\left(\sum_{k=1}^{n-1}\Es(\xi_k|\F_{k-1})\right).$$
 
 Now, almost surely, either $<M>_\iy<\iy$ or $M_n/<M>_n\des_{n\to\iy}0$, by Lévy's strong law of large numbers (see for instance \cite{williams}, Chapter 12, § 14), with yields the result.
 \ep
\setlength{\textheight}{21cm}

\medskip

{\bf Acknowledgment.} I would like to thank Terry Lyons for encouraging me to give a graduate course on self-interacting processes in Oxford in 2005 and 2007, and Nathanael Enriquez and Christophe Sabot for their invitation to lecture a minicourse on reinforced walks in Aussois in June 2009, both of which initiated these notes. 

\footnotesize
\bibliographystyle{plain}
\bibliography{localization1}

\end{document}